\documentclass[12pt]{article}
\usepackage{amsfonts}

\usepackage{latexsym}
\usepackage{amssymb}
\usepackage{amsmath}
\usepackage{amsthm}
\usepackage[mathscr]{eucal}
\usepackage{graphicx}
\usepackage{hyperref}
\usepackage{caption}
\usepackage{subcaption}

\newtheorem{Lemma}{Lemma}
\newtheorem{Proposition}{Proposition}
\newtheorem{Theorem}[Lemma]{Theorem}

\newtheorem{Definition}{Definition}

\renewcommand{\qed}{\hfill{\ \ \rule{2mm}{2mm}} \vspace{0.2in}}

\newcommand{\ind}{1\hspace{-2.3mm}{1}}

\begin{document}

\title{Euclidean minimum spanning trees with location dependent and power weighted edges}
%\titlerunning{Euclidean minimum spanning trees}
%\authorrunning{G. Ganesan}
\author{ \textbf{Ghurumuruhan Ganesan}
\thanks{E-Mail: \texttt{gganesan82@gmail.com} } \\
%EndAName
\ \\
Institute of Mathematical Sciences, HBNI, Chennai.}
\date{}
\maketitle

\begin{abstract}
Consider~\(n\) nodes~\(\{X_i\}_{1 \leq i \leq n}\)
independently distributed in the unit square~\(S,\)
each according to a distribution~\(f\)
and let~\(K_n\) be the complete graph formed by joining
each pair of nodes by a straight line segment.
For every edge~\(e\) in~\(K_n\) we associate a weight~\(w(e)\)
that may depend on the \emph{individual locations} of the endvertices
of~\(e.\) Denoting~\(MST_n\) to be the minimum weight of a spanning tree of~\(K_n\)
and assuming an equivalence condition on the weight function~\(w(.),\)
we prove that~\(MST_n\) appropriately
scaled and centred converges to zero a.s.\
and in mean as~\(n \rightarrow \infty.\)
We also obtain upper and lower bound deviation estimates for~\(MST_n.\)

%\(\frac{1}{\sqrt{n}}(MST_n - \mathbb{E}MST_n)\)

 %and analogous results hold for minimum spanning trees and paths.

%In this paper, we study the structure of left-right crossings
%of the random geometric graph \(G = G(n,r_n)\) of \(n\) nodes
%uniformly distributed in \(S = [0,1]^2\) with \(r_n = \epsilon\sqrt{\frac{\log{n}}{n}}\)
%for some \(\epsilon > 0.\) Tiling \(S\) horizontally and
%vertically into rectangles of length \(1\) and width \(Mr_n,\) we
%show that each rectangle has a left-right crossing of edges with
%high probability if \(M\) is sufficiently large.
%We call the resulting subgraph to be a ``backbone" of \(G.\)

%The techniques we use to construct the backbone has quite a few applications.
%As a first, we show that the diameter of second largest component in \(G\)
%is \(O(1)\) with high probability. Secondly,
\vspace{0.1in} \noindent \textbf{Key words:} Minimum spanning tree, location dependent edge weights, edge weight exponent.

\vspace{0.1in} \noindent \textbf{AMS 2000 Subject Classification:} Primary:
60J10, 60K35; Secondary: 60C05, 62E10, 90B15, 91D30.
\end{abstract}

\bigskip

\setcounter{equation}{0}
\renewcommand\theequation{\thesection.\arabic{equation}}
\section{Introduction}\label{intro}
The study of the minimum weight spanning trees of a graph is of great practical importance
and many algorithms have been proposed over the years
for various kinds of graphs. For example, the well-known Kruskal's algorithm (Cormen et al (2009))
iteratively adds edges to a sequence of increasing subtree of the original graph
until a spanning tree is obtained with the constraint that no cycle is created
in any of the iterations. The spanning tree with minimum weight so obtained
is usually called the \emph{Minimum Spanning Tree} (MST).

We are interested in MSTs of Euclidean random graphs whose nodes
are randomly distributed in the unit square and whose edges are assigned weights
related to the Euclidean length. When the weight of an edge equals its Euclidean length raised to a positive power,
we refer to the resulting MSTs as power weighted Euclidean MSTs.

One of the main objects of interest in the study of power weighted Euclidean MST is its total weight: How does it scale with
the number of nodes and the power weight exponent and what are its convergence properties?
Analytical results for such MSTs  have been studied extensively before (see Steele (1988, 1993), Kesten and Lee (1996),
Penrose and Yukich (2003) and references therein).
For example, Steele (1988) uses edge counting techniques to obtain variance estimates for the MST weight
and Kesten and Lee (1996) use martingale methods to obtain central limit theorems (CLTs)
for the MST weight, appropriately scaled and centred.
Penrose and Yukich (2003) use coupling arguments to obtain weak laws
for functionals of point processes thereby including the MST as a special case.
Recently Chatterjee and Sen (2017) used percolation theoretic arguments to study
convergence rate of the CLTs for Euclidean MSTs.

In this paper, we study MSTs of Euclidean random graphs whose edge weights depend on the \emph{individual locations} of the endvertices.
Such MSTs frequently arise in practice and as an example, suppose it is required to establish a fully connected wireless communication network among the set of nodes randomly located in a geographical area. An edge between two nodes signifies a potential communication link with the corresponding weight being proportional to the cost of physically installing the link. If some areas within the area are highly populated, it might actually be cost effective to set up links in such ``hotspots". Thus the edge weights are location dependent and it is of interest to estimate the minimum cost of setting up such a network. Throughout we use the unit square for illustrating our results and the analysis holds analogous for other regular shapes like rectangle or circle.

%Our proof techniques is also used to analyse MST of nodes distributed throughout the unit square with location dependent edge weights and we obtain
%variance estimates and almost sure convergence results for the length of the MST, suitably scaled and centred %(see Theorem~\ref{var_mst_thm}).

%Using our proof techniques, we also obtain variance estimates

In the rest of this section, we describe the model under consideration and describe an important property regarding the maximum vertex degree, that distinguishes location dependent MSTs from their homogenous (i.e. location independent) counterparts. We then describe the outline and the main objectives of the paper towards understanding the properties of location dependent MSTs.

\subsection*{\em Model Description}
Let~\(f\) be any distribution on the unit square~\(S\) satisfying such that
\begin{equation}\label{f_eq}
\epsilon_1 \leq \inf_{x \in S} f(x) \leq \sup_{x \in S} f(x) \leq \epsilon_2
\end{equation}
for some constants~\(0 < \epsilon_1 \leq \epsilon_2 < \infty.\)
Throughout all constants are independent of~\(n.\)

%where~\(vol(A)\) refers to the area of~\(A.\)

%For \(1 \leq i \leq n,\) we define \(X_i\) on the probability space \((S, {\cal B}(S), \mathbb{P}_i),\) where \(\mathbb{P}_i\) denotes the Lebesgue measure on \(S.\)

Let~\(\{X_i\}_{i \geq 1}\) be independently
and identically distributed (i.i.d.) with the distribution~\(f(.)\) defined on the probability space \((\Omega, {\cal F}, \mathbb{P}).\)
For~\(n \geq 1,\) let~\(K_n = K(X_1,\ldots,X_n)\) be the complete graph whose edges are
obtained by connecting each pair of nodes~\(X_i\) and~\(X_j\) by the straight line segment~\((X_i,X_j)\)
with~\(X_i\) and~\(X_j\) as endvertices.

A path~\({\cal P} = (Y_1,\ldots,Y_t)\) containing~\(t\) distinct vertices is a subgraph of~\(K_n\)
with vertex set~\(\{Y_{j}\}_{1 \leq j \leq t}\)
and edge set~\(\{(Y_{j},Y_{{j+1}})\}_{1 \leq j \leq t-1}.\)
The nodes~\(Y_1\) and~\(Y_t\) are said to be \emph{connected} by
edges of the path~\({\cal P}.\)

A subgraph~\({\cal T}\) of~\(K_n\)
with vertex set~\(\{Y_i\}_{1 \leq i \leq t}\) and edge set~\(E_{\cal T}\)
is said to be a \emph{tree} (Steele (1988)) if
the following two conditions hold:\\
\((1)\) The graph~\({\cal T}\) is connected; i.e., any two nodes
in~\({\cal T}\) are connected by a path containing only edges in~\(E_{\cal T}.\)\\
\((2)\) The graph~\({\cal T}\) is acyclic; i.e., no subgraph of~\({\cal T}\) is a cycle.\\
The tree~\({\cal T}\) is said to be a \emph{spanning tree} of~\(K_n\) if~\({\cal T}\)
contains all the~\(n\) nodes\\\(\{X_k\}_{1 \leq k \leq n}.\)

In what follows we assign weights to edges of the graph~\(K_n\) and study minimum weight spanning trees.

\subsection*{\em Minimum spanning trees}
For points~\(x,y \in S,\) we let~\(d(x,y)\) denote the Euclidean distance between~\(x\) and~\(y\) and let~\(h : S \times S \rightarrow (0,\infty)\) be a deterministic measurable function satisfying
\begin{equation}\label{eq_met1}
c_1 d(x,y) \leq h(x,y) =h(y,x) \leq c_2 d(x,y) %CHNG NOTATION PRKVMM+etC..
\end{equation}
for some positive constants~\(c_1,c_2.\) For a constant~\(\alpha > 0\) and for~\(1 \leq i < j \leq n\) we let~\(h^{\alpha}(X_i,X_j) = h^{\alpha}(e)\) denote the \emph{weight} of the edge~\(e = (X_i,X_j)\) with corresponding edge weight exponent~\(\alpha.\) We simply refer to~\(h^{\alpha}(e)\) as the weight of the edge~\(e\) and define~\(d(e) = d(X_i,X_j)\) to be the Euclidean length of~\(e.\) Unless mentioned otherwise, all weights are with respect to the edge weight exponent~\(\alpha\) and all lengths are Euclidean. %all weights are with respect to the exponent~\(\alpha\) and

For a tree~\({\cal T}\) with vertex set~\(\{Y_1,\ldots,Y_t\},\) the weight of~\({\cal T}\) is the sum of the weights of the edges in~\({\cal T};\) i.e.,~\(W({\cal T}) := \sum_{e \in {\cal T}} h^{\alpha}(e).\) Let~\({\cal T}_{n}\) be a spanning tree of~\(K_n\) containing all the nodes~\(\{X_i\}_{1 \leq i \leq n}\) and satisfying
\begin{equation}\label{min_weight_tree}
MST_n = W({\cal T}_{n}) := \min_{{\cal T}} W({\cal T}),
\end{equation}
where the minimum is taken over all spanning trees~\({\cal T}.\) We denote~\({\cal T}_{n}\) to be the \emph{minimum spanning tree} (MST) with corresponding weight~\(MST_n.\) If there is more than one choice for~\({\cal T}_{n},\) we choose one according to a deterministic rule.

%By definition we always have that~\(MST_n \leq MSP_n \leq TSP_n\) since any spanning cycle contains a spanning path which in turn is a spanning tree.

%WRT PROPO...here that degree could be arbitrarily large...

Let~\(d(X_i) = d\left(X_i,{\cal T}_n,h\right)\) denote the degree of the node~\(X_i,1 \leq i \leq n\) in the MST~\({\cal T}_n\) with edge weight function~\(h\) so that~\(\sum_{i=1}^{n}d(X_i) = 2(n~-~1),\) since~\({\cal T}_n\) contains~\(n-1\) edges. Taking expectations, we therefore get that\\\(n\mathbb{E}d(X_1) = 2(n-1)\) and so~\(\mathbb{E}d(X_1) \leq 2.\) In case the edge weight function~\(h~=~d,\) the Euclidean length, then the weight of any edge does not depend on the location of its endvertices and the \emph{maximum} degree of any node in~\({\cal T}_n\) is at most~\(6,\) almost surely. This is because the angle between any two edges in the MST sharing an endvertex, cannot be more than 60 degrees (see Aldous and Steele (1992)). This geometric property frequently occurs in the study of the properties of the random MST (Kesten and  Lee~(1996), Yukich (2000)).

For edge weight functions~\(h\) that are location dependent, it is possible for the maximum degree of a vertex in the MST to take arbitrary large values with positive probability as stated in the following result.
\begin{Proposition}\label{prop1} For every integer~\(K \geq 2\) there exists an edge weight function~\(h = h_K\) satisfying~(\ref{eq_met1}), a sequence~\(n_l \rightarrow \infty\) as~\(l \rightarrow \infty\) and a constant~\(\epsilon_0  = \epsilon_0(K) > 0\) such that
\begin{equation}\label{deg_x1}
\mathbb{P}\left(\max_{1 \leq i \leq n_l} d\left(X_i,{\cal T}_{n_l},h\right) \geq K\right) \geq \epsilon_0
\end{equation}
for all~\(l \geq 1.\)
\end{Proposition}
In light of Proposition~\ref{prop1}, we are therefore interested to study the effect of large vertex degrees on the total weight of location dependent MSTs. We demonstrate via common properties like variance, deviation estimates and almost sure convergence, that the total weight behaves roughly the same way as in the location independent case. %in spite of arbitrarily large vertex degrees.

\subsection*{\em Paper outline}
The paper is organized as follows. In Section~\ref{pf_prop}, we prove Proposition~\ref{prop1} and in Section~\ref{mst_dev}, we state and prove deviation estimates for the weight of location dependent MSTs
and also obtain upper and lower bounds for their expected value. We use stochastic domination and coupling techniques
to obtain the deviation estimates and provide numerical simulations
to illustrate the bounds obtained. We show that the MST weight~\(MST_{n}\) defined in~(\ref{min_weight_tree})
still remains of the order of~\(n^{1-\frac{\alpha}{2}}\) as in the location independent case.

Next in Section~\ref{mst_var_up}, we use the deviation estimates obtained in Section~\ref{mst_dev} to find upper bounds on the variance of
location dependent MSTs. The main tools we use in this section are one node difference estimates and the martingale difference method
and we also briefly explain using Propostion~\ref{prop1}, how the analysis varies between the location dependent and independent cases.
As before the variance upper bound is of the order of~\(n^{1-\alpha}\) and the same order as in the location independent case.

In Section~\ref{mst_var_low}, we obtain lower bounds on the variance of MST weight
by studying occurrence of predetermined nice configurations that differ in the location
of at most one node. In particular combining with the variance upper bound obtained in Section~\ref{mst_var_low},
we obtain that the variance \emph{is} of the order of~\(n^{1-\alpha}.\)

Finally in Sections~\ref{mst_conv} and~\ref{mst_uni}, for completeness,
we study a.s.\ convergence for location dependent MSTs and uniform MSTs using subsequence arguments
and the deviation and variance estimates obtained in the previous sections.

\setcounter{equation}{0}
\renewcommand\theequation{\thesection.\arabic{equation}}
\section{Proof of Proposition~\ref{prop1}}\label{pf_prop}
\begin{figure}
\centering
\begin{subfigure}{0.5\textwidth}
\centering
\includegraphics[width=2.5in, trim= 50 250 50 100, clip=true]{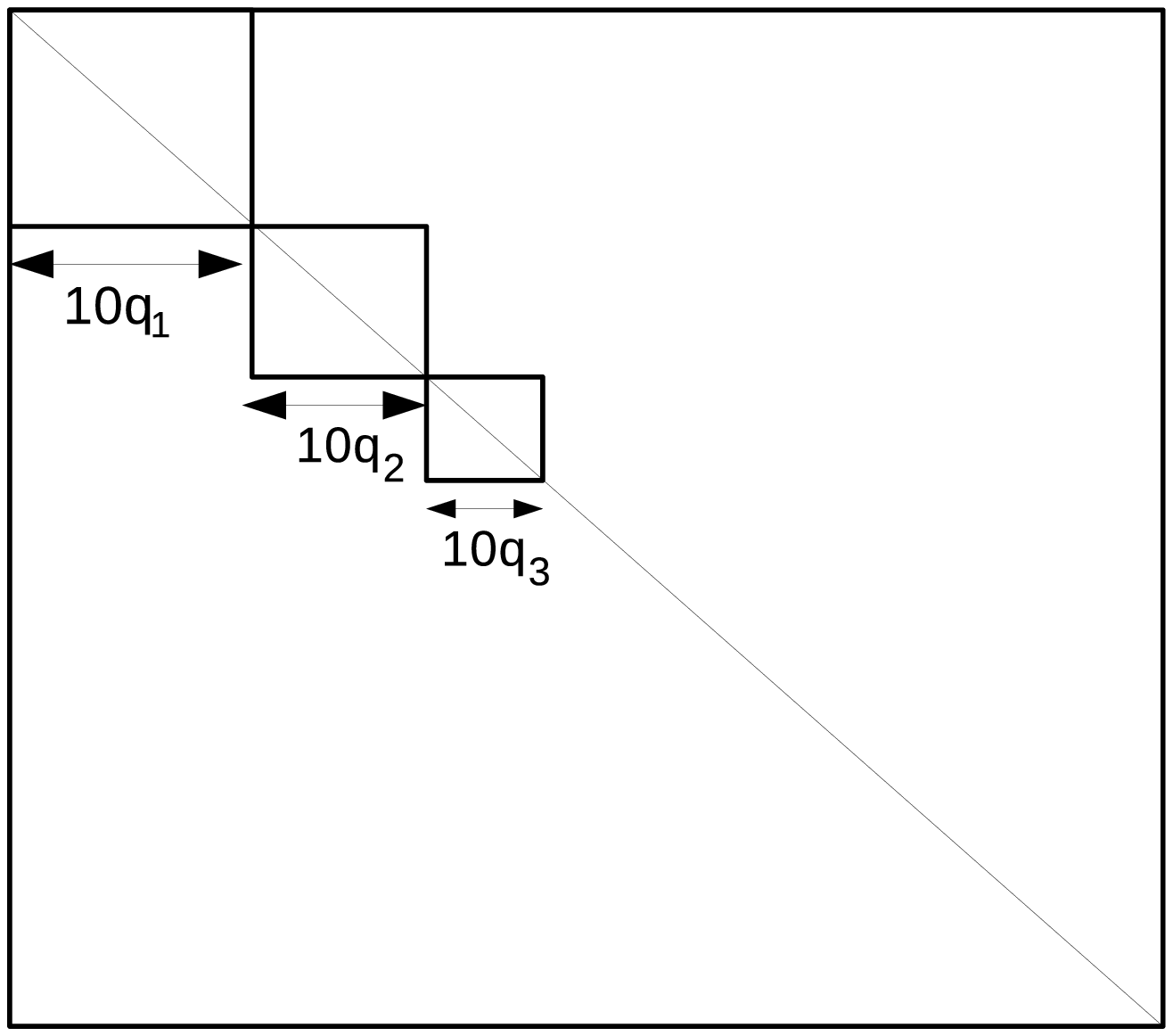}
  \caption{} %\caption{Tiling~\(S\) into~\(r_n \times r_n\) squares with an inter-square distance of~\(s_n.\)}%\nonumber%{t1}
\end{subfigure}% seems this is important....
\begin{subfigure}{0.5\textwidth}
\centering
   \includegraphics[width=2.5in, trim= 0 250 50 100, clip=true]{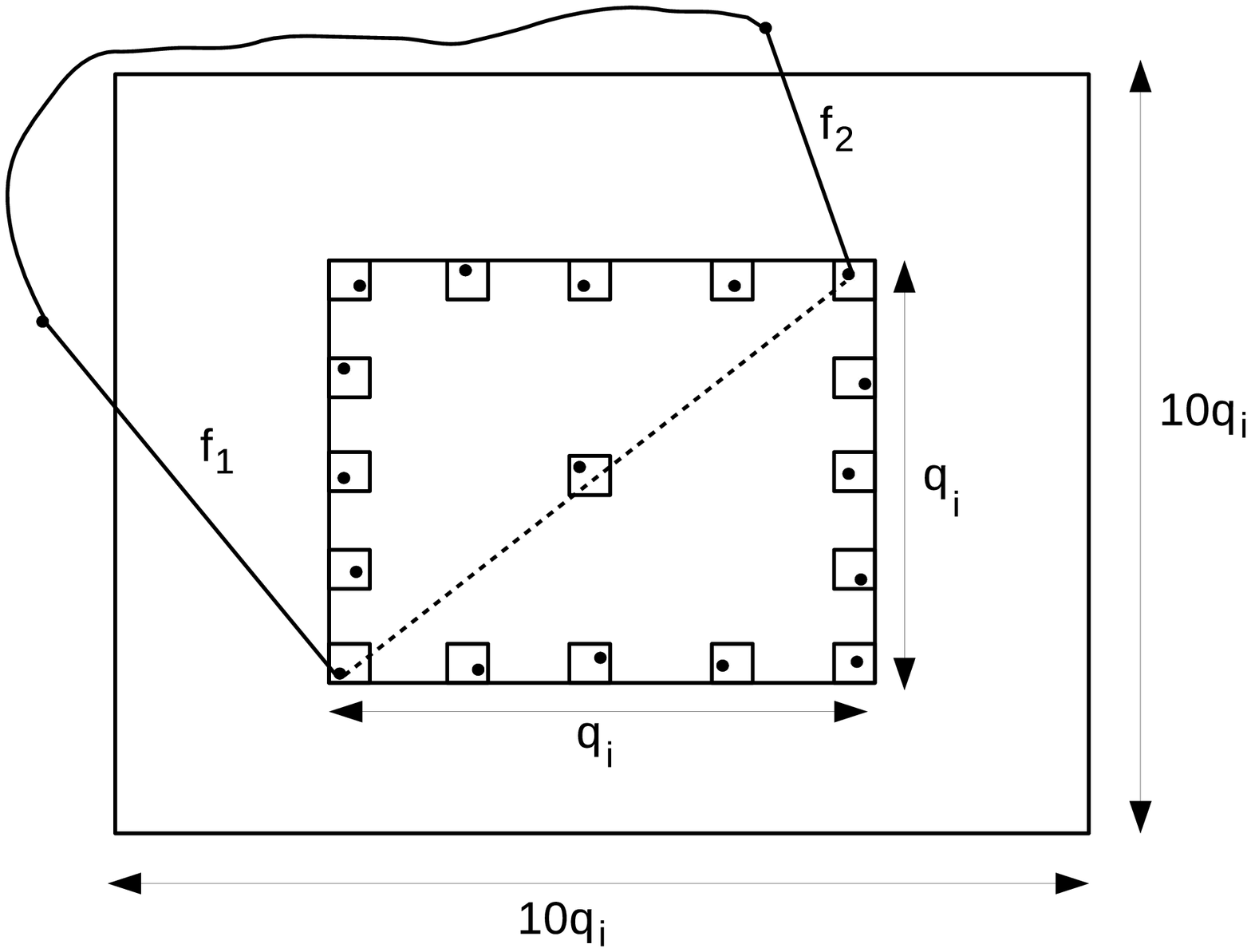}
   \caption{} %\caption{The induced subgraph~\(G({\cal D}_n(N))\) corresponding to the centres of the grey squares in~\((a).\)} %\nonumber%{t2}
  \end{subfigure}

\caption{\((a)\) The~\(10 q_i \times 10q_i\) squares placed along the diagonal of the unit square. \((b)\) The~\(q_i \times q_i\) square~\(S_i\) along with the smaller subsquares.}
\label{diag_squares}
\end{figure}

We start with some preliminary computations. For~\(i \geq 1\) let~\(n_i = D\cdot i^3\) for some constant~\(D > 0\) to be determined later and let~\(q_i = \frac{2K-1}{\sqrt{n_i}}.\) Place disjoint~\(10q_i \times 10q_i\) squares~\(S^{big}_i\) along the diagonal of the unit square~\(S\) as shown in Figure~\ref{diag_squares}~\((a)\) and choose~\(D > 0\) large enough so that the total sum length of the diagonals of the squares in~\(\{S^{big}_i\}\) is at most the length of the diagonal of~\(S;\) i.e., we choose~\(D\) large enough so that \[\sum_{i \geq 1}10 q_i \sqrt{2} \leq  \sqrt{2}.\]

Let~\(S_i\) be the~\(q_i \times q_i\) subsquare with the same centre as~\(S^{big}_i\) as shown in Figure~\ref{diag_squares}~\((b)\) and consider~\(K\) small subsquares each of size~\(\frac{1}{\sqrt{n_i}} \times \frac{1}{\sqrt{n_i}}\) placed~\(x_i = \frac{1}{\sqrt{n_i}}\) apart on each of the four sides of~\(S_i.\) The total number of the small subsquares is~\(4(K-1),\) which we label as~\(S_i(1),\ldots,S_i(4K-4).\) The central small square (labeled~\(S_i(0)\)) has the same centre as the square~\(S_i\) and is also of size~\(\frac{1}{\sqrt{n_i}} \times \frac{1}{\sqrt{n_i}}.\)

Let~\(c_1\) and~\(c_2\) be positive constants such that~\(c_1 < \frac{c_2}{8K}\) and define the edge weight function~\(h\) as:\\
\begin{equation} \label{h_def22}
h(u,v) = \left\{
\begin{array}{cc}
c_1 d(u,v) & \text{ if } u \in \cup_j \{S_j(0)\} \text{ or } v \in \cup_j \{S_j(0)\} \\
 & \\
c_2 d(u,v) & \text{ otherwise},
\end{array}
\right.
\end{equation}
where~\(d(u,v)\) represents the Euclidean distance between~\(u\) and~\(v\) as before. By choice of~\(h\) in~(\ref{h_def22}), the edge weight containing a node is less if the node is present in one of the central squares and
larger otherwise.

%~\(h(u,v) = c_1d(u,v)\) if either~\(u\) or~\(v\) is in~\(\cup_j \{S_j(0)\}\) and~\(h(u,v)=c_2d(u,v),\) otherwise.

For~\(i \geq 1,\) let~\({\cal T}_{n_i}\) be the MST of the nodes~\(\{X_k\}_{1 \leq k \leq n_i}\) as defined in~(\ref{min_weight_tree}) with edge weight function~\(h\) as in~(\ref{h_def22}). To identify nodes in~\({\cal T}_{n_i}\) with large degree we define~\(E_K(i)\)
to be the event that there is exactly one node~\(v_i \in \{X_k\}_{1 \leq k \leq n_i}\) present in~\(S_i(l)\) for each~\(0 \leq l \leq 4K-4\) and the rest of~\(S_{i}^{big}\) is empty. We recall from the first paragraph of this section that~\(S_{i}^{big}\) is the\\\(10q_i \times 10q_i\) square with the same centre as~\(S_i(0).\)

The following Lemma implies that the event~\(E_K(i)\) occurs with positive probability and that if~\(E_K(i)\) occurs, then the node~\(v_0\) present in the central square~\(S_i(0)\) has large degree.
\begin{Lemma}\label{lem_prop1} There exists a constant~\(\epsilon_0> 0\) depending only on~\(K,\epsilon_1\) and~\(\epsilon_2\) such that
\begin{equation}\label{pek_est}
\mathbb{P}(E_K(i)) \geq \epsilon_0.
\end{equation}

If~\(E_K(i)\) occurs and if~\({\cal T}_{loc}\) denotes the induced subgraph of~\({\cal T}_{n_i}\) formed by the nodes in~\(\{v_j\}_{0 \leq j \leq 4K-4},\) then~\({\cal T}_{loc}\) is a tree (more specifically, a star graph) with edge set~\(\{(v_0,v_j)\}_{1 \leq j \leq 4K-4}.\)
\end{Lemma}
From Lemma~\ref{lem_prop1}, we therefore get that the degree of~\(v_0\) in the MST~\({\cal T}_{n_i}\) is~\(4K-4\) if the event~\(E_K(i)\) occurs. From the probability estimate~(\ref{pek_est}), we therefore get Proposition~\ref{prop1}.

We now prove Lemma~\ref{lem_prop1} beginning with~(\ref{pek_est}).\\
\emph{Proof of~(\ref{pek_est}) in Lemma~\ref{lem_prop1}}: Any node of~\(\{X_k\}_{1 \leq k \leq n_i}\) is present within the square~\(S_i(l)\) with probability~\(p_i(l):= \int_{S_i(l)}f(x)dx\) and so
\begin{equation}\label{pek}
\mathbb{P}(E_K(i)) = {n_i \choose 4K-3} (4K-3)!\prod_{l=0}^{4K-4} p_i(l) \left(\int_{S \setminus S_{i}^{big}}f(x)dx\right)^{n_i-4K+3}.
\end{equation}
The area of the square~\(S_i(l)\) equals~\(\frac{1}{n_i}\) and so~\(p_i(l) \geq \frac{\epsilon_1}{n_i}\) using the bounds in~(\ref{f_eq}). Similarly,~\(\int_{S_{i}^{big}}f(x)dx \leq \epsilon_2 100q_i^2 = \frac{C}{n_i}\) where~\(C = 100\epsilon_2(2K-1)^2\) and so
\begin{equation}\label{pek2}
\mathbb{P}(E_K(i)) \geq {n_i \choose 4K-3}(4K-3)!  \left(\frac{\epsilon_1}{n_i}\right)^{4K-3} \left(1-\frac{C}{n_i}\right)^{n_i}.
\end{equation}

For~\(a > 2b\) we have~\({a \choose b}b! \geq (a-b)^{b} \geq \frac{a^b}{2^{b}}.\) Therefore choosing~\(i\) larger if necessary so that~\(n_i \geq 8K-6\) and setting~\(a = n_i\) and~\(b = 4K-3,\) we get from~(\ref{pek2}) that
\begin{equation}\label{pek3}
\mathbb{P}(E_K(i)) \geq \left(\frac{\epsilon_1}{2}\right)^{4K-3}\cdot \left(1-\frac{C}{n_i}\right)^{n_i}.
\end{equation}
To evaluate the final term in~(\ref{pek3}), we use~\(1-y \geq e^{-2y}\) for all~\(0 < y < \frac{1}{2}.\) To see this estimate is true
we write~\(\log(1-y) =-y - \sum_{k \geq 2} \frac{y^{k}}{k}\) and use~\(y < \frac{1}{2}\)
to get that~\[\sum_{k \geq 2} \frac{y^{k}}{k} \leq \frac{1}{2}\sum_{k \geq 2} y^{k} = \frac{y^2}{2(1-y)} \leq y^2  < y.\]
Choosing~\(i\) large enough so that~\(y = \frac{(4K-3)\epsilon_2}{n_i} < \frac{1}{2},\)
we then get from~(\ref{pek3}) that
\begin{equation} \nonumber
\mathbb{P}(E_K(i)) \geq \left(\frac{\epsilon_1}{2}\right)^{4K-3}\cdot e^{-2C},
\end{equation}
proving~(\ref{pek_est}).~\(\qed\)

\emph{Proof of rest of Lemma~\ref{lem_prop1}}: Let~\(S_i^{big}\) be the~\(10q_i \times 10q_i\) square in Figure~\ref{diag_squares}\((b).\) First we  show that~\({\cal T}_{loc}\) is a tree by a contradiction argument. Suppose for example that two nodes in opposite~\(\frac{1}{\sqrt{n_i}} \times \frac{1}{\sqrt{n_i}}\) squares are joined by a path~\(P\) in~\({\cal T}_{n_i}\) containing at least one vertex not in~\(\{v_i\}\) as in Figure~\ref{diag_squares}\((b).\) Since~\(E_K(i)\) occurs, the rest of~\(S_i^{big}\) not containing~\(\cup_l S_i(l)\) is empty and so the path~\(P\) contains at least two ``long" edges~\(f_1 = (u_1,x_1),f_2 = (u_2,x_2)\) with endvertices outside~\(S_i^{big}\) as shown in Figure~\ref{diag_squares}\((b).\) The nodes~\(u_1\) and~\(u_2\) both belong to~\(\{v_i\}_{1 \leq i \leq 4K-4}\) and the dotted edge~\((u_1,u_2)\) is not present in the MST~\({\cal T}_{n_i}\) because otherwise, we would get a cycle.

The weights of the edges~\(f_1\) and~\((u_1,u_2)\) equal~\((c_2d(u_1,x_1))^{\alpha}\) and\\\((c_2d(u_1,u_2))^{\alpha},\) respectively, by definition of the edge weight function in~(\ref{h_def22}). The length~\(d(u_1,x_1)\) of the edge~\(f_1 = (u_1,x_1)\) is at least~\(\frac{9q_i}{2},\)  the width of the annulus~\(S_i^{big} \setminus S_i\) and the length of~\((u_1,u_2)\) at most~\(q_i \sqrt{2}.\) This implies that removing the edge~\(f_1\) and adding~\((u_1,u_2),\) we would get an MST with weight strictly larger than~\({\cal T}_{n_i},\) a contradiction. This proves that~\({\cal T}_{loc}\) is a tree.

Next, to see that~\({\cal T}_{loc}\) is a star graph, we use the definition of the edge weight function~(\ref{h_def22}) and obtain that the weight of any edge of the form~\((v_0,v_j)\) is at most
\begin{equation}\label{huy_one}
(c_1 \cdot d(v_0,v_j))^{\alpha} \leq \left(c_1 \cdot \left(\frac{q_i}{\sqrt{2}} +\frac{\sqrt{2}}{\sqrt{n_i}}\right)  \right)^{\alpha} \leq \left(\frac{c_1\cdot 4K }{\sqrt{n_i}}\right)^{\alpha}
\end{equation}
since~\(q_i = \frac{2K-1}{\sqrt{n_i}}\) (see first paragraph of this Section). Similarly, the weight of any edge of the form~\((v_{s},v_t), s, t \neq 0\) is at least~\(\left(\frac{c_2}{\sqrt{n_i}}\right)^{\alpha}.\) Since~\(c_1 < \frac{c_2}{8K}\) (see~(\ref{h_def22})), we get that~\((v_s,v_t)\) has larger weight that~\((v_0,v_j).\) Thus the minimum weight tree containing all the nodes~\(\{v_j\}_{0 \leq j \leq 4K-4}\) is the  star graph with vertex set~\(\{(v_0,v_j)\}_{0 \leq j \leq 4K-4}.\)~\(\qed\)

\setcounter{equation}{0}
\renewcommand\theequation{\thesection.\arabic{equation}}
\section{Deviation estimates for the MST}\label{mst_dev}
As a first step in the study of properties of location dependent MSTs, we obtain deviation estimates in this Section.
The bounds obtained are of same order and we also use these estimates in later Sections for the study of variance and almost sure convergence.

We begin with a couple of preliminary definitions. Let~\(\epsilon_1,\epsilon_2\) be as in~(\ref{f_eq}) and set~\(\delta = \delta(\alpha) = \epsilon_1\) if the edge weight exponent~\(\alpha \leq 1\) and~\(\delta = \epsilon_2\) if~\(\alpha > 1.\) Recalling that~\(c_1\) and~\(c_2\) are the bounds for the edge weights as in~(\ref{eq_met1}), we define for~\(A > 0\) the terms~\(C_1(A)=C_1(A,c_1,c_2,\epsilon_1,\epsilon_2,\alpha)\) and~\(C_2(A) = C_2(A,c_1,c_2,\epsilon_1,\epsilon_2,\alpha)\) as
\begin{eqnarray}
C_1(A)  &:=& \frac{(c_1A)^{\alpha}}{2A^2}(1-e^{-\epsilon_1 A^2})e^{-8\epsilon_2 A^2} \text{ and } \nonumber\\
C_2(A)  &:=& (2c_2A)^{\alpha}\left(1 + \frac{\mathbb{E}T^{\alpha}}{A^2}\right), \label{c12def}
\end{eqnarray}
where~\(T\) is a geometric random variable with success parameter~\(p = 1-e^{-\delta A^2},\) independent of the node locations~\(\{X_i\};\)
i.e.,~\(\mathbb{P}(T = k) = (1-p)^{k-1}p\) for all integers~\(k \geq 1.\)  Letting~\(MST_n\)
be the MST weight as defined in~(\ref{min_weight_tree}), we have the following main result.
\begin{Theorem}\label{thm_mst_dev} Let~\(\alpha > 0\) be the edge weight exponent. For every~\(A > 0\) and integer~\(k \geq 1\) and all~\(n \geq n_0(A,k)\) large,
\begin{equation}\label{mst_low_bounds}
\mathbb{P}\left(MST_n \geq C_1(A) n^{1-\frac{\alpha}{2}}\left(1-\frac{4\sqrt{A}}{n^{1/4}}\right)\right) \geq 1-e^{-n^{1/3}},
\end{equation}
\begin{equation}\label{mst_up_bounds}
\mathbb{P}\left(MST_n \leq C_2(A) n^{1-\frac{\alpha}{2}}\left(1+ \frac{2}{n^{1/16}}\right)\right) \geq 1-\frac{1}{n^{2k}}
\end{equation}
and
\begin{equation}
C^{k}_1(A) \left(1-\frac{36k\sqrt{A}}{n^{1/4}}\right)\leq \mathbb{E}\left(\frac{MST^{k}_n}{n^{k\left(1-\frac{\alpha}{2}\right)}}\right) \leq C^{k}_2(A)\left(1+ \frac{2k}{n^{1/16}}\right). \label{exp_mst_bound}
\end{equation}
\end{Theorem}

\subsection*{\em Remarks on Theorem~\ref{thm_mst_dev}}
From Theorem~\ref{thm_mst_dev}, we see that the weights of the MST in the location dependent case, is of the same order~\(n^{1-\frac{\alpha}{2}}\) as in the location independent case (see Steele (1988)). In fact, using~(\ref{exp_mst_bound}) we get that the normalized MST weight~\(\frac{\mathbb{E} MST_n}{n^{1-\frac{\alpha}{2}}}\) in fact satisfies
\begin{equation}\label{gr}
c_1^{\alpha} \cdot \beta_{low}(\alpha) \leq \liminf_n \frac{\mathbb{E} MST_n}{n^{1-\frac{\alpha}{2}}} \leq \limsup_n \frac{\mathbb{E} MST_n}{n^{1-\frac{\alpha}{2}}} \leq c_2^{\alpha}\cdot \beta_{up}(\alpha),
\end{equation}
where
\begin{equation}\label{beta_low_def}
\beta_{low}(\alpha) = \beta_{low}(\alpha, \epsilon_1,\epsilon_2) := \sup_{A > 0} \frac{A^{\alpha}}{2A^2}(1-e^{-\epsilon_1 A^2})e^{-8\epsilon_2 A^2},
\end{equation}
\begin{equation}\label{beta_up_def}
\beta_{up}(\alpha) = \beta_{up}(\alpha,\epsilon_1,\epsilon_2) := \inf_{A > 0} (2A)^{\alpha}\left(1 + \frac{\mathbb{E}{T}^{\alpha}}{A^2}\right),
\end{equation}
and~\(T\) is a geometric random variable with success parameter~\(p = 1-e^{-\delta A^2}\) (see~(\ref{c12def})).

For the case of homogenous distribution~\(\epsilon_1 = \epsilon_2 = 1\) we get
\begin{equation}\label{d_one}
\beta_{low}(\alpha) = \frac{1}{2}\sup_{A > 0}A^{\alpha-2}(1-e^{-A^2})e^{-8A^2} > 0
\end{equation}
and
\begin{equation}\label{d_two}
\beta_{up}(\alpha) = \inf_{A > 0} (2A)^{\alpha} \left(1+\frac{\mathbb{E}T^{\alpha}}{A^2}\right) < \infty,
\end{equation}
where~\(T\) is a geometric random variable with success parameter~\(p = 1-e^{-A^2}.\)

\begin{figure}[tbp]
\centering
%\fbox{
\includegraphics[width=3.5in, clip=true]{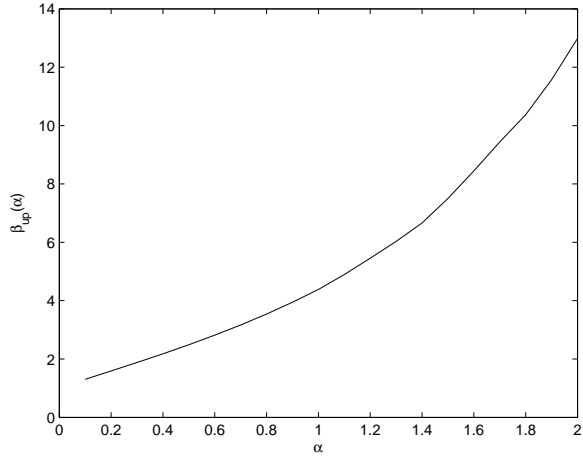}
%}
\caption{Plot of~\(\beta_{up}(\alpha)\) as a function of~\(\alpha\) for the homogenous case~\(\epsilon_1 = \epsilon_2 = 1.\)}
\label{fig_up}
\end{figure}

\begin{figure}[tbp]
\centering
%\fbox{
\includegraphics[width=3.5in, clip=true]{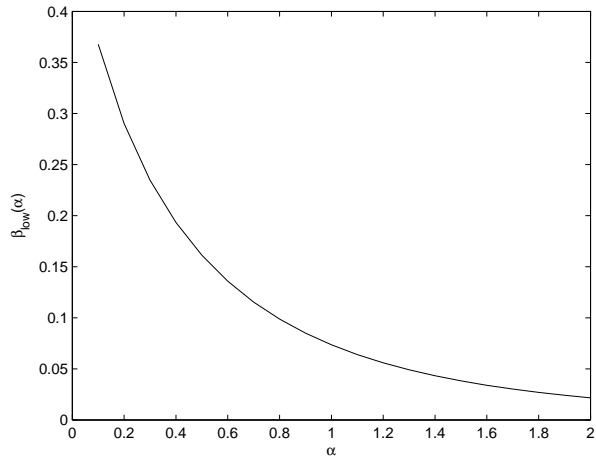}
%}
\caption{Plot of~\(\beta_{low}(\alpha)\) as a function of~\(\alpha\) for the homogenous case~\(\epsilon_1 = \epsilon_2 = 1.\)}
\label{fig_down}
\end{figure}

For illustration, we plot~\(\beta_{low}(\alpha)\) and~\(\beta_{up}(\alpha)\) as a function of~\(\alpha\) in Figures~\ref{fig_down} and~\ref{fig_up}, respectively. As we see from the figures~\(\beta_{up}(\alpha)\) increases with~\(\alpha\) and~\(\beta_{low}(\alpha)\) decreases with~\(\alpha.\) Expressions~(\ref{d_one}) and~(\ref{d_two}) also allow us to numerically evaluate the bounds in~(\ref{gr}) for various values of~\(\alpha.\) For example, for~\(\alpha = 1,\) we get that
\[\beta_{low}(1) = \frac{1}{2} \sup_{A > 0} A^{-1} (1-e^{-A^2})e^{-8A^2} \approx 0.0735633 > \frac{1}{20}\] and since~\(\mathbb{E}T = \frac{1}{p} = \frac{1}{1-e^{-A^2}},\) we get that
\[\beta_{up}(1) = \inf_{A > 0} 2A\left(1+ \frac{1}{A^2(1-e^{-A^2})}\right) \approx 4.46256 < 5.\] Substituting these bounds back in~(\ref{gr})
we find that
\begin{equation}\label{gr2}
\frac{c_1^{\alpha}}{20} \leq \liminf_n \frac{\mathbb{E} MST_n}{n^{1-\frac{\alpha}{2}}} \leq \limsup_n \frac{\mathbb{E} MST_n}{n^{1-\frac{\alpha}{2}}} \leq 5 \cdot c_2^{\alpha}
\end{equation}
and this obtains bounds for the scaled MST weights in terms of the edge weight function parameters~\(c_1\) and~\(c_2.\)

Similarly for~\(\alpha = 2,\) we get that~
\[\beta_{low}(2) = \frac{1}{2} \sup_{A > 0} (1-e^{-A^2})e^{-8A^2} \approx 0.0216525\] and since~\(\mathbb{E}T^2 = \frac{2-p}{p^2} = \frac{1+e^{-A^2}}{(1-e^{-A^2})^2},\) we get that
\[\beta_{up}(2) = \inf_{A > 0} (2A)^2\left(1+ \frac{1+e^{-A^2}}{A^2(1-e^{-A^2})^2}\right)   \approx 13.8772\]
and we have analogous bounds as in~(\ref{gr2}) in this case as well.

In general,~\(\beta_{low}(\alpha)\) and~\(\beta_{up}(\alpha)\) in~(\ref{gr}) can be evaluated to get some knowledge of the dependence of the normalized MST weight on the node distribution parameters~\(\epsilon_1,\epsilon_2\) (see~(\ref{f_eq})) and the edge weight function parameters~\(c_1,c_2\) (see~(\ref{eq_met1})). For example, suppose the distribution~\(f = \frac{1}{2}\) on~\([0,0.5]^2\) and~\(f = \frac{7}{6}\) on the remaining area in~\(S\) so that~\(\int f = 1.\) In this case~\(\epsilon_1 = \frac{1}{2}\) and~\(\epsilon_2 = \frac{7}{6}\) and we directly evaluate~(\ref{beta_low_def}) and~(\ref{beta_up_def}) to get that~\(\beta_{low}(1) \approx 0.0346363\) and~\(\beta_{up}(1) \approx 4.92912.\)

As a final remark, we provide a simple upper bound for~\(\mathbb{E}T^{\alpha}\) in order to obtain quick evaluations of~\(\beta_{up}(\alpha)\) in~(\ref{beta_up_def}). Using~\(\mathbb{P}(T \geq k) = (1-p)^{k-1} \leq e^{-p(k-1)}\) where~\(p = 1-e^{-\delta A^2},\) we see that relation~(\ref{x_dist}) in Appendix is satisfied and so letting~\(r\) be the smallest integer greater than or equal~\(\alpha\) we have from~(\ref{disc_tel}) that
\[\mathbb{E}T^{\alpha} \leq \mathbb{E}T^{r} \leq \frac{r!}{(1-e^{-p})^{r}}.\] Plugging this estimate into~(\ref{beta_up_def}) provides
an upper bound for~\(\beta_{up}(\alpha).\)

\subsection*{\em Proof outline and Poissonization}
The rough idea for the proof of Theorem~\ref{thm_mst_dev} is as follows. We first tile the unit square into small subsquares and for obtaining the lower bound, we determine configurations within these subsquares that result in long edges. For the upper bound, we explicitly construct a spanning tree whose weight is at most a constant multiple of the weight of the MST, with high probability, i.e., with probability one as~\(n \rightarrow \infty.\)

To prove the deviation estimates in Theorem~\ref{thm_mst_dev}, we use Poissonization and let~\({\cal P}\) be a Poisson process
in the unit square~\(S\) with intensity~\(nf(.).\) We join  each pair
of nodes by a straight line segment and denote the resulting complete graph as~\(K^{(P)}_n.\)
As in~(\ref{min_weight_tree}), we let~\(MST_n^{(P)}\) be the weight of the MST of~\(K^{(P)}_n\) and first find deviation estimates for~\(MST_{n}^{(P)}.\) We then use dePoissonization to obtain corresponding deviation estimates for~\(MST_n,\) the weight of the MST of the complete graph~\(K_n\) as defined in~(\ref{min_weight_tree}).

\begin{figure}[tbp]
\centering
%\fbox{
\includegraphics[width=3in, trim= 20 200 50 110, clip=true]{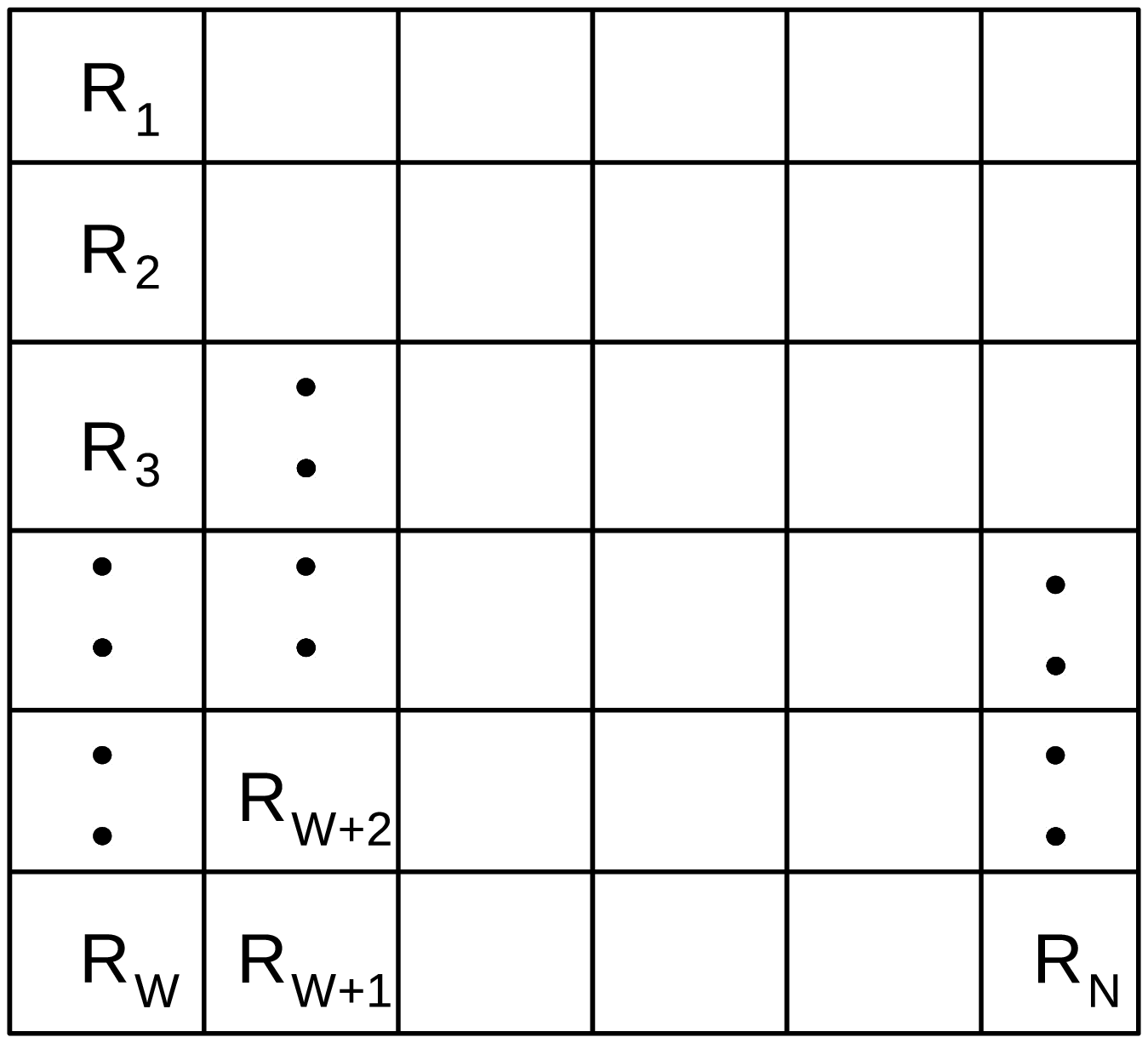}
%}
\caption{Tiling the unit square into~\(N = \frac{n}{A^2}\) smaller~\(\frac{A}{\sqrt{n}} \times \frac{A}{\sqrt{n}}\) squares~\(\{R_l\}_{1 \leq l \leq \frac{n}{A^2}}.\)}
\label{fig_squares}
\end{figure}

For a  real number~\(A > 0,\) we tile the unit square~\(S\) into small~\(\frac{A(n)}{\sqrt{n}} \times \frac{A(n)}{\sqrt{n}}\) squares~\(\{R_i\}_{1 \leq i \leq \frac{n}{A^2(n)}}\) where~\(A(n) \in \left[A, A + \frac{1}{\log{n}}\right]\) is chosen such that~\(\frac{\sqrt{n}}{A(n)}\) is an integer. This is possible since
\begin{equation}
\frac{\sqrt{n}}{A} - \frac{\sqrt{n}}{A + (\log{n})^{-1}} = \frac{\sqrt{n}}{\log{n}}\cdot \frac{1}{A(A+(\log{n})^{-1})} \geq \frac{\sqrt{n}}{2A^2\log{n}} \label{poss}
\end{equation}
for all~\(n\) large. For notational simplicity, we denote~\(A(n)\) as~\(A\) henceforth and label the squares as in Figure~\ref{fig_squares} so that~\(R_i\) and~\(R_{i+1}\) share an edge for~\(1 \leq i \leq \frac{n}{A^2}-1.\)

\subsection*{\em Lower deviation bounds~(\ref{mst_low_bounds})}
For~\(1 \leq i \leq \frac{n}{A^2}\) let~\(E(R_i)\) denote the event that the~\(\frac{A}{\sqrt{n}} \times \frac{A}{\sqrt{n}}\) square~\(R_i\) is occupied i.e., contains at least one node of~\({\cal P},\) and all squares sharing a corner with~\(R_i\) are empty. If~\(E(R_i)\) occurs, then there is at least one edge in the MST of~\(K^{(P)}_n\) with one endvertex in~\(R_i\) and other endvertex in a square not sharing a corner with~\(R_i.\) Such an edge has a Euclidean length of at least~\(\frac{A}{\sqrt{n}}\) and so a weight of at least~\(\left(\frac{c_1A}{\sqrt{n}}\right)^{\alpha}\) (see~(\ref{eq_met1})). Consequently
\begin{equation}\label{low_bound}
MST^{(P)}_n \geq \frac{1}{2} \cdot \sum_{i=1}^{\frac{n}{A^2}} \left(\frac{c_1A}{\sqrt{n}}\right)^{\alpha} \ind(E(R_i)) = \frac{1}{2} \cdot\left(\frac{c_1A}{\sqrt{n}}\right)^{\alpha}\cdot G_{\alpha},
\end{equation}
where~\(G_{\alpha} := \sum_{i=1}^{\frac{n}{A^2}} \ind(E(R_i))\) and the factor~\(\frac{1}{2}\) occurs, since each edge is counted twice in the summation.

To estimate~\(G_{\alpha},\) we would like to split it into sums of independent r.v.s using the following construction. For a square~\(R_i,\) let~\({\cal N}(R_i)\) be the set of all squares sharing a corner with~\(R_i,\) including~\(R_i.\) If~\(R_i\) does not intersect the sides of the unit square~\(S,\) then there are~\(9\) squares in~\({\cal N}(R_i)\) and if~\(R_j\) is another square such that~\({\cal N}(R_i) \cap {\cal N}(R_j) = \emptyset,\) then the corresponding events~\(E(R_i)\) and~\(E(R_j)\) are independent, by Poisson property. We therefore extract nine disjoint subsets~\(\{{\cal U}_l\}_{1 \leq l \leq 9}\) of~\(\{R_i\}\) with the following properties:\\
\((A)\) If~\(R_i,R_j \in {\cal U}_l,\) then~\(\#{\cal N}(R_i) = \#{\cal N}(R_j) = 9\) and~\({\cal N}(R_i) \cap {\cal N}(R_j) = \emptyset.\)\\
\((B)\) The number of squares~\(\#{\cal U}_l \geq \frac{n}{9A^2} - \frac{4\sqrt{n}}{A}\) for each~\(1 \leq l \leq 9.\)\\
This is possible since there are at most~\(\frac{4\sqrt{n}}{A} - 4 < \frac{4\sqrt{n}}{A}\) squares in~\(\{R_k\}\) intersecting the sides of the unit square~\(S\) and the total number of squares in~\(\{R_k\}\) is~\(\frac{n}{A^2}.\)

We now write~\(G_{\alpha} = \sum_{i=1}^{\frac{n}{A^2}} \ind(E(R_i)) \geq \sum_{l=1}^{9} \sum_{R_i \in {\cal U}_l} \ind(E(R_i)),\) where each inner summation on the right side is a sum of independent Bernoulli random variables, which we bound via standard deviation estimates. Indeed for~\(1 \leq l \leq 9\) and~\(R_i \in {\cal U}_l,\) the number of nodes~\(N(R_i)\) is Poisson distributed with mean~\(n\int_{R_i}f(x)dx \in [\epsilon_1 A^2, \epsilon_2 A^2]\) (see~(\ref{f_eq})) and so~\(R_i\) is occupied with probability at least~\(1-e^{-\epsilon_1 A^2}.\) Also each of the eight squares sharing a corner with~\(R_i\) is empty with probability at least~\(e^{-\epsilon_2 A^2},\) implying that~\(\mathbb{P}(E(R_i)) \geq (1-e^{-\epsilon_1 A^2})e^{-8\epsilon_2 A^2}.\) Using the standard deviation estimate~(\ref{std_dev_down}) of Lemma~\ref{app_lem} in Appendix with~\(\mu_1 = (1-e^{-\epsilon_1 A^2})e^{-8\epsilon_2 A^2}, m = \frac{n}{9A^2}-\frac{4\sqrt{n}}{A}\) and~\(\epsilon = \frac{1}{m^{1/4}}\) we then get that
\begin{equation}\label{ul_est}
\mathbb{P}_0\left(\sum_{R_i \in {\cal U}_l} \ind(E(R_i)) \geq (1-\epsilon)\left(\frac{n}{9A^2}-\frac{4\sqrt{n}}{A}\right)(1-e^{-\epsilon_1 A^2})e^{-8\epsilon_2 A^2}\right) \geq 1-e^{-D_1 \epsilon^2 n}
\end{equation}
for some constant~\(D_1 > 0\) not depending on~\(l.\) Since~\(m^{1/4} < \left(\frac{n}{9A^2}\right)^{1/4},\) we get that~\(D_1\epsilon^2 n \geq 2D_2 \sqrt{n}\)  for some constant~\(D_2 > 0\) and since~\(m^{1/4} > \left(\frac{n}{10A^2}\right)^{1/4}\) for all~\(n\) large, we have \[(1-\epsilon)\left(\frac{n}{9A^2} - \frac{4\sqrt{n}}{A}\right) \geq \frac{n}{9A^2} - \frac{4\sqrt{n}}{A} - \frac{n}{A^2 m^{1/4}} \geq \frac{n}{9A^2}\left(1 - \frac{36\sqrt{A}}{n^{1/4}}\right)\] for all~\(n\) large.

Letting
\[E_{low} := \left\{G_{\alpha}\geq (1-e^{-\epsilon_1 A^2})e^{-8\epsilon_2 A^2}\frac{n}{A^2}\left(1 - \frac{36\sqrt{A}}{n^{1/4}}\right)\right\},\] we get from~(\ref{ul_est}) that~\(\mathbb{P}_0(E_{low}) \geq 1-9e^{-2D_2 \sqrt{n}}\) and moreover, from~(\ref{low_bound}) we also get that \[MST^{(P)}_n \ind(E_{low}) \geq \Delta_n := C_1(A) n^{1-\frac{\alpha}{2}} \left(1-\frac{36\sqrt{A}}{n^{1/4}}\right),\] where~\(C_1(A)\) is as defined in~(\ref{c12def}). From the estimate for the probability of the event~\(E_{low}\) above we therefore get~\(\mathbb{P}_0\left(MST_n^{(P)}\geq \Delta_n\right) \geq 1- 9e^{-2D_2 \sqrt{n}}\) for all~\(n\) large. To convert the estimate from Poisson to the Binomial process, we let~\(E_P := \left\{MST_n^{(P)} \geq \Delta_n\right\}, E := \left\{MST_n \geq \Delta_n\right\}\) and use the dePoissonization formula (Bradonjic et al. (2010))
\begin{equation}\label{de_poiss_ax}
\mathbb{P}(E) \geq 1- D \sqrt{n} \mathbb{P}(E^c_P)
\end{equation}
for some constant~\(D > 0\) to get that~\(\mathbb{P}(E) \geq 1- D\sqrt{n}e^{-2D_2\sqrt{n}} \geq 1~-~e^{-D_2 \sqrt{n}}\) for all~\(n\) large. This proves~(\ref{mst_low_bounds}) and so using~\[\mathbb{E}MST^{k}_n \geq \mathbb{E}MST^{k}_n \ind\left(MST_n \geq \Delta_n\right) \geq \Delta^{k}_n \left(1-e^{-D_2\sqrt{n}}\right)\] and
\[\left(1-\frac{36\sqrt{A}}{n^{1/4}}\right)^{k} (1-e^{-D_2 \sqrt{n}}) \geq \left(1-\frac{36k\sqrt{A}}{n^{1/4}}\right)(1-e^{-D_2 \sqrt{n}})  \geq 1-\frac{37k\sqrt{A}}{n^{1/4}}\] for all~\(n\) large, we also obtain the lower bound on the expectation in~(\ref{exp_mst_bound}).

To see that~(\ref{de_poiss_ax}) is true, we let~\(N_P\) denote the random number of nodes of~\({\cal P}\) in all the squares~\(\{S_j\}\) so that~\(\mathbb{E}_0 N_P = n\) and~\(\mathbb{P}_0(N_P=n) = e^{-n}\frac{n^{n}}{n!} \geq \frac{D_1}{\sqrt{n}}\) for
some constant~\(D_1 > 0,\) using the Stirling formula. Given~\(N_P = n,\) the nodes of~\({\cal P}\)
are i.i.d.\ with distribution~\(f(.)\) as defined in~(\ref{f_eq}); i.e., \(\mathbb{P}_0(E_P^c|N_P = n)  = \mathbb{P}(E^c)\) and so
\[\mathbb{P}_0(E_P^c) \geq \mathbb{P}_0(E_P^c|N_P = n) \mathbb{P}_0(N_P = n) =
\mathbb{P}(E^c) \mathbb{P}_0(N_P = n) \geq \mathbb{P}(E^c)\frac{D_1}{\sqrt{n}},\] proving~(\ref{de_poiss_ax}).~\(\qed\)

%We also need another estimate that merges trees within small squares in~\(S.\)

\subsection*{\em Upper deviation bounds~(\ref{mst_up_bounds})}
As before, we first obtain upper bounds for the Poissonized MST~\(MST^{(P)}_n.\) The idea is to first connect all the nodes within each square~\(R_j\) to get a collection of subtrees and then join all these subtrees together to get an overall spanning tree. Suppose there is at least one node of the Poisson process~\({\cal P}\) in the unit square~\(S\) and let~\(R_{i_1},R_{i_2},\ldots,R_{i_Q}, 1 \leq i_1 < i_2 <\ldots < i_Q \leq \frac{n}{A^2}, Q \leq \frac{n}{A^2}\) be the~\(\frac{A}{\sqrt{n}} \times \frac{A}{\sqrt{n}}\) squares containing all the nodes of~\({\cal P}.\) For~\(1 \leq j \leq Q,\) let~\({\cal T}_{i_j}\) be any spanning tree containing all the nodes of~\(R_{i_j}\) and for~\(1 \leq j \leq Q-1\) let~\(e_{j+1}\) be any edge joining some node in~\(R_{i_j}\) and some node in~\(R_{i_{j+1}}\) so that
the union~\({\cal T}_{uni} := \cup_{1 \leq j \leq Q} {\cal T}_{i_j} \cup \cup_{2 \leq l \leq Q} \{e_l\}\) is a spanning tree of the complete graph~\(K^{(P)}_n.\) The weight~\(W({\cal T}_{uni}) \geq MST_n^{(P)}\) and so it suffices to upper bound~\(W({\cal T}_{uni}).\)

For~\(1 \leq j \leq Q,\) there are~\(N(R_{i_j})\) nodes of the Poisson process in the square~\(R_{i_j}\) and any two such nodes are connected by an edge of length at most~\(\frac{A\sqrt{2}}{\sqrt{n}}.\) Therefore the spanning tree~\({\cal T}_{i_j}\) has a total weight of at most\\\(N(R_{i_j}) \cdot \left(\frac{c_2A\sqrt{2}}{\sqrt{n}}\right)^{\alpha},\) using the bounds for the function~\(h\) in~(\ref{eq_met1}). The edge~\(e_{j+1}\) that connects some node in~\(R_{i_j}\) with some node of~\(R_{i_{j+1}}\) has a Euclidean length of at most~\(\frac{2T_{j+1}A}{\sqrt{n}}\) where~\(T_{j+1} := i_{j+1} - i_{j}\) and therefore has a weight of at most~\(\left(\frac{2c_2T_{j+1}A}{\sqrt{n}}\right)^{\alpha},\) by~(\ref{eq_met1}). In effect,
\begin{eqnarray}
W\left({\cal T}_{uni}\right) &\leq& \sum_{j=1}^{Q} N(R_{i_j}) \cdot \left(\frac{c_2A\sqrt{2}}{\sqrt{n}}\right)^{\alpha} +
\sum_{j=1}^{Q-1} \left(\frac{2c_2T_{j+1}A}{\sqrt{n}}\right)^{\alpha} \nonumber\\
&=& \sum_{i=1}^{\frac{n}{A^2}} N(R_{i}) \cdot \left(\frac{c_2A\sqrt{2}}{\sqrt{n}}\right)^{\alpha} +
\sum_{j=2}^{Q} \left(\frac{2c_2T_{j}A}{\sqrt{n}}\right)^{\alpha}. \nonumber
\end{eqnarray}
Setting~\(T_1 := i_1-1,T_{Q+1} := \frac{n}{A^2}-i_Q\) and~\(S_{\alpha} := \sum_{j=1}^{Q+1} T_j^{\alpha}\) we then get
\begin{eqnarray}\label{up_bd1}
MST_n^{(P)} \leq W\left({\cal T}_{uni}\right) \leq \left(\frac{2c_2A}{\sqrt{n}}\right)^{\alpha} \left(\sum_{i=1}^{\frac{n}{A^2}} N(R_{i}) + S_{\alpha}\right).
\end{eqnarray}

The first sum~\(\sum_{i=1}^{\frac{n}{A^2}} N(R_{i})\) in the right side of~(\ref{up_bd1}) is a Poisson random variable with mean~\(n\) since this denotes the total number of nodes of the Poisson process in the unit square. Using the deviation estimate~(\ref{std_dev_up}) in Appendix with~\(m=1,\mu_2 = n\) and~\(\epsilon = \frac{\log{n}}{\sqrt{n}},\) we have that~
\begin{equation}\label{e_node_est}
\mathbb{P}_0\left(\sum_{i=1}^{\frac{n}{A^2}} N(R_{i}) \leq n\left(1+\frac{\log{n}}{\sqrt{n}}\right)\right) \geq 1-e^{-D(\log{n})^2}
\end{equation}
for some constant~\(D > 0.\) Setting~\(E_{node} :=  \left\{\sum_{i=1}^{\frac{n}{A^2}} N(R_i) \leq n\left(1+\frac{\log{n}}{\sqrt{n}}\right)\right\},\) we get from~(\ref{up_bd1}) that
\begin{equation}\label{up_bd2}
MST_n^{(P)} \ind(E_{node})  \leq \left(\frac{2c_2A}{\sqrt{n}}\right)^{\alpha} \left(n\left(1+\frac{\log{n}}{\sqrt{n}}\right)  + S_{\alpha}\right).
\end{equation}

The second term~\(S_{\alpha}\) in~(\ref{up_bd2}) is well-defined for any configuration \(\omega\) of the Poisson process provided we set~\(S_{\alpha}(\omega_0) = \left(\frac{n}{A^2}-1\right)^{\alpha}\) where~\(\omega_0\) is the configuration containing no node of the Poisson process in the unit square~\(S.\) The following Lemma obtains an estimate on~\(S_{\alpha}.\)
\begin{Lemma} \label{s_alpha_lemma} Let~\(T\) be a geometric random variable with success parameter\\\(p = 1-e^{-A^2\delta}\) independent of the node locations~\(\{X_i\}.\) For every even integer~\(m \geq 1\) and every~\(A > 0\) and for all~\(n \geq n_0(m,A)\) large, we have
\begin{equation}
\mathbb{P}_0\left(S_{\alpha} \leq  \left(1+\frac{1}{n^{1/16}}\right)\frac{n}{A^2} \mathbb{E}T^{\alpha}\right)\geq 1-\frac{D}{n^{7m/16}}, \label{til_est2}
\end{equation}
where~\(D > 0\) is a constant.
\end{Lemma}
Since~\(S_{\alpha}\) is not an i.i.d.\ sum, we use coupling techniques to estimate~\(S_{\alpha}\) and prove Lemma~\ref{s_alpha_lemma} at the end of this section.

We continue with the proof of the deviation upper bound. Let~\(m \geq 2\) be an even integer to be determined later. Letting~\(E_{S}\) be the event on the left side of~(\ref{til_est2}) we have from the previously computed upper bound on~\(MST_n^{(P)}\) (see~(\ref{up_bd2})) that
\begin{eqnarray} \nonumber
MST_n^{(P)} \ind(E_{node} \cap E_{S})   &\leq& \left(\frac{2c_2A}{\sqrt{n}}\right)^{\alpha} \left(n\left(1+\frac{\log{n}}{\sqrt{n}}\right)   + \left(1+\frac{1}{n^{1/16}}\right)\frac{n}{A^2}\mathbb{E}T^{\alpha}\right) \nonumber\\
&\leq& C_2(A) n^{1-\frac{\alpha}{2}}\left(1+ \frac{1}{n^{1/16}}\right) \label{up_bd3}
\end{eqnarray}
where~\(C_2(A)\) is as in~(\ref{c12def}) and the final estimate in~(\ref{up_bd3}) is obtained using~\(\frac{2\log{n}}{\sqrt{n}} \leq \frac{1}{n^{1/16}}\) for all~\(n\) large.  From~(\ref{up_bd3}) and the estimates for the events\\
\(E_{node}\) and~\(E_{S}\) from~(\ref{e_node_est}) and~(\ref{til_est2}), respectively, we have
\begin{eqnarray}
\mathbb{P}_0\left(MST_n^{(P)} \leq C_2(A) n^{1-\frac{\alpha}{2}} \left(1+\frac{1}{n^{1/16}}\right)\right) &\geq& 1 - e^{-D(\log{n})^2} -  \frac{D}{n^{7m/16}} \nonumber\\
&\geq& 1-  \frac{2D}{n^{7m/16}} \label{main_poi_up}
\end{eqnarray}
for all~\(n\) large and some constant~\(D > 0.\)

From~(\ref{main_poi_up}) and the dePoissonization formula~(\ref{de_poiss_ax}), we obtain
\begin{equation}
\mathbb{P}\left(MST_n \leq C_2(A) n^{1-\frac{\alpha}{2}} \left(1+\frac{1}{n^{1/16}}\right)\right) \geq 1- \frac{2D \sqrt{n}}{n^{7m/16}}. \label{temp11}
\end{equation}
Choosing~\(m\) large enough such that~\(\frac{7m}{16}-\frac{1}{2} \geq k,\) we obtain the estimate in~(\ref{mst_up_bounds}).

For bounding the expectation of~\(MST_n^{k},\) we let~\(\Delta_n = C_2(A) n^{1-\frac{\alpha}{2}} \left(1+\frac{1}{n^{1/16}}\right)\) and write
\begin{eqnarray}
\mathbb{E}MST^{k}_n &=& \mathbb{E}MST^{k}_n\ind(MST_n \leq \Delta_n) + \mathbb{E}MST^{k}_n \ind(MST_n > \Delta_n) \nonumber\\
&\leq& \Delta^{k}_n + \mathbb{E}MST^{k}_n \ind(MST_n > \Delta_n) \label{term2}
\end{eqnarray}
To estimate~\(\Delta_n^{k},\) we use~\(\left(1+\frac{1}{n^{1/16}}\right)^{k} \leq e^{k/n^{1/16}} \leq 1+\frac{2k}{n^{1/16}}\) for all~\(n\) large and get that~\(\Delta_n^{k} \leq C_2^{k}(A)n^{k\left(1-\frac{\alpha}{2}\right)}\left(1+\frac{2k}{n^{1/16}}\right)\) for all~\(n\) large. For the second term in~(\ref{term2}), we use the estimate~\(MST_n \leq n \cdot \left(c_2\sqrt{2}\right)^{\alpha}\) since there are~\(n\) edges in the spanning tree, each such edge has an Euclidean length of at most~\(\sqrt{2}\) and so the weight of any edge is at most~\((c_2 \sqrt{2})^{\alpha},\) using~(\ref{eq_met1}). Letting~\(\theta_m = \frac{7m}{16}-\frac{1}{2} - \frac{\alpha k}{2}\) and using the probability estimate~(\ref{temp11}) we then get that
\[\mathbb{E} MST^{k}_n \leq \Delta^{k}_n + \frac{3D_2n^{k}\sqrt{n}}{n^{7m/16}} \leq C^{k}_2(A)n^{k\left(1-\frac{\alpha}{2}\right)}\left(1+\frac{2k}{n^{1/16}} + \frac{D_3}{n^{\theta_{m}}}\right) \] for all~\(n\) large and some constant~\(D_3 > 0.\) Choosing~\(m\) larger if necessary so that~\(\theta_m \geq 1 > \frac{1}{16},\) we obtain the expectation upper bound in~(\ref{exp_mst_bound}).~\(\qed\)

\emph{Proof of Lemma~\ref{s_alpha_lemma}}: We show that~\(S_{\alpha}(\omega)\) is monotonic in~\(\omega\) in the sense that adding more nodes increases~\(S_{\alpha}\) if~\(\alpha \leq 1\) and decreases~\(S_{\alpha}\) if~\(\alpha > 1.\) This then allows us to use coupling and upper bound~\(S_{\alpha}\) by simply considering homogenous Poisson processes.

\underline{\emph{Monotonicity of~\(S_{\alpha}\)}}: We recall that~\(\omega_0\) is the configuration containing no node of the Poisson process in the unit square. For a configuration~\(\omega \neq \omega_0\) let~\(1 \leq i_1(\omega) < \ldots < i_Q(\omega) \leq \frac{n}{A^2}\) be the indices of the squares in~\(\{R_j\}\) containing at least one node of the Poisson process~\({\cal P}.\) Letting~\(i_0(\omega) = 1\) and~\(i_{Q+1}(\omega) = \frac{n}{A^2}\) we have~\(S_{\alpha}(\omega) = \sum_{j=0}^{Q}(i_{j+1}(\omega)-i_j(\omega))^{\alpha}.\) Suppose~\(\omega' = \omega \cup \{x\}\) is obtained by adding a single extra node at~\(x \in R_{j_0}\) for some~\(1 \leq j_0 \leq \frac{n}{A^2}.\) If~\(j_0 \in \{i_k(\omega)\}_{0 \leq k \leq Q+1},\) then~\(S_{\alpha}(\omega') = S_{\alpha}(\omega).\) Else there exists~\(0 \leq a \leq Q\) such that~\(i_a(\omega)  < j_0 < i_{a+1}(\omega)\) and so~\[S_{\alpha}(\omega') = S_{\alpha}(\omega) + (i_{a+1}(\omega)-j_0)^{\alpha} +(j_0-i_a(\omega))^{\alpha} - (i_{a+1}(\omega)-i_a(\omega))^{\alpha}.\]

If~\(\alpha \leq 1\) then using~\(a^{\alpha} + b^{\alpha} \geq (a+b)^{\alpha} \) for positive numbers~\(a,b\) we get that~\(S_{\alpha}(\omega') \geq S_{\alpha}(\omega).\) If~\(\alpha > 1\) then~\(a^{\alpha}  + b^{\alpha} \leq (a+b)^{\alpha}\) and so~\(S_{\alpha}(\omega') \leq S_{\alpha}(\omega).\) This monotonicity property together with coupling allows us to upper bound~\(S_{\alpha}\) as follows. Letting~\(\delta = \epsilon_2\) if~\(\alpha \leq 1\) and~\(\delta = \epsilon_1\) if~\(\alpha > 1,\) we let~\({\cal P}_{\delta}\) be a homogenous Poisson process of intensity~\(\delta n\) on the unit square~\(S,\) defined on the probability space~\((\Omega_{\delta},{\cal F}_{\delta},\mathbb{P}_{\delta}).\)

Let~\(F_{\delta}\) denote the event that there is at least node of~\({\cal P}_{\delta}\) in the unit square~\(S\) and set~\(S^{(\delta)}_{\alpha} := \left(\frac{n}{A^2}-1\right)^{\alpha}\) if~\(F_{\delta}\) does not occur. If~\(F_{\delta}\) occurs, then as before let~\(\{i^{(\delta)}_{j}\}_{1 \leq j \leq Q_{\delta}}\) be the indices of the squares in~\(\{R_j\}\) containing at least one node of~\({\cal P}_{\delta}.\) Moreover, let~\(T^{(\delta)}_{j+1} := i^{(\delta)}_{j+1} - i^{(\delta)}_j\) for~\(1 \leq j \leq Q_{\delta}\) and set~\(T^{(\delta)}_1 := i^{(\delta)}_1-1\) and~\(T^{(\delta)}_{Q_{\delta}+1} := \frac{n}{A^2}-i^{(\delta)}_{Q_{\delta}}.\) Defining~\(S^{(\delta)}_{\alpha} := \sum_{j=1}^{Q_{\delta}+1} \left(T^{(\delta)}_j\right)^{\alpha}\) in this case,  we have for any~\(x > 0\) that
\begin{equation}\label{mon_salpha}
\mathbb{P}_{\delta}\left(S^{(\delta)}_{\alpha} < x\right) \leq \mathbb{P}_0\left(S_{\alpha} < x\right).
\end{equation}
The proof of~(\ref{mon_salpha}) follows from standard coupling arguments and for completeness, we provide a proof in Appendix.

%for any~\(x > 0.\) %and some constant~\( D> 0\) not depending on~\(n\) or~\(x.\)

To estimate~\(S^{(\delta)}_{\alpha}\) we let~\(N^{(\delta)}(R_i),1 \leq i \leq \frac{n}{A^2},\) be the random number of nodes of~\({\cal P}_{\delta}\) in the square~\(R_i.\) The random variables~\(\{N^{(\delta)}(R_i)\}\) are i.i.d.\ Poisson distributed each with mean~\(A^2\delta.\)  For~\(i \geq \frac{n}{A^2}+1,\) we define~\(N^{(\delta)}(R_i)\) to be i.i.d.\ Poisson random variables with mean~\(A^2\delta,\) that are also independent of~\(\{N^{(\delta)}(R_i)\}_{1 \leq i \leq \frac{n}{A^2}}.\) Without loss of generality, we associate the probability measure~\(\mathbb{P}_{\delta}\) for the random variables~\(\{N^{(\delta)}(R_i)\}_{i \geq \frac{n}{A^2} +1}\) as well.

Let~\(\tilde{T}_1 := \min\{j \geq 1 : N^{(\delta)}(R_j) \geq 1\}\)  and for~\(j \geq 2,\) let~\[\tilde{T}_j := \min\{k \geq \tilde{T}_{j-1}+1 : N^{(\delta)}(R_k) \geq 1\}-\tilde{T}_{j-1}.\] The random variables~\(\{\tilde{T}_i\}\) are nearly the same as~\(\{T^{(\delta)}_i\}\) in the following sense: Suppose the event~\(F_{\delta}\) occurs so that there is at least one node of~\({\cal P}_{\delta}\) in the unit square. This means that~\(1 \leq Q_{\delta} \leq \frac{n}{A^2}\) and so~\(T^{(\delta)}_1 = i_1-1 =\tilde{T}_1-1, T^{(\delta)}_j = \tilde{T}_j\) for~\(2 \leq j \leq Q_{\delta}\) and~\(T^{(\delta)}_{Q_{\delta}+1} \leq \tilde{T}_{Q_{\delta}+1}.\) Consequently
\begin{equation}\label{dom1}
S_{\alpha}^{(\delta)} \ind(F_{\delta}) \leq \sum_{i=1}^{Q_{\delta}+1} \tilde{T}^{\alpha}_i\ind(F_{\delta}) \leq \sum_{i=1}^{\frac{n}{A^2}+1} \tilde{T}^{\alpha}_i\ind(F_{\delta}) \leq \sum_{i=1}^{\frac{n}{A^2}+1} \tilde{T}^{\alpha}_i,
\end{equation}
since~\(Q_{\delta} \leq \frac{n}{A^2}.\)

From~(\ref{dom1}), it suffices to estimate the sum on the right side to find an upper bound for~\(S^{(\delta)}_{\alpha}.\) The advantage of~(\ref{dom1}) is that~\(\{\tilde{T}_i\}\) are i.i.d.\ geometric random variables with success parameter~\(p = 1-e^{-A^2\delta}\) and so all moments of~\(\tilde{T}_i^{\alpha}\) exist. Letting~\(\beta_i = \left(\tilde{T}_i^{\alpha} - \mathbb{E}_{\delta}\tilde{T}_i^{\alpha}\right)\) and~\(\beta_{tot} = \sum_{i=1}^{\frac{n}{A^2}+1} \beta_i \) we obtain for an even integer constant~\(m\) that
\begin{equation}\label{e_del}
\mathbb{E}_{\delta}\beta_{tot}^m  = \mathbb{E}_{\delta} \sum_{(i_1,\ldots,i_m)} \beta_{i_1}\ldots\beta_{i_m}.
\end{equation}
For a tuple~\((i_1,\ldots,i_m)\) let~\(\{j_1,\ldots,j_w\}\) be the distinct integers in~\(\{i_1,\ldots,i_m\}\) with corresponding multiplicities~\(l_1,\ldots,l_w\) so that~\[\mathbb{E}_{\delta} \beta_{i_1}\ldots\beta_{i_m} = \mathbb{E}_{\delta} \beta_{j_1}^{l_1} \ldots \beta_{j_w}^{l_w} = \prod_{k=1}^{w} \mathbb{E}_{\delta} \beta_{j_k}^{l_k}.\] If~\(l_k = 1\) for some~\(1 \leq k \leq w,\) then~\(\mathbb{E}_{\delta} \beta_{i_1}\ldots\beta_{i_m} = 0\) and so for any non zero term in the summation in~(\ref{e_del}), there are at most~\(\frac{m}{2}\) distinct terms in~\(\{i_1,\ldots,i_m\}.\) This implies that
\[\mathbb{E}_{\delta}\beta_{tot}^{m} \leq D(m) {n \choose m/2} \leq D(m) n^{m/2}\] for some constant~\(D(m) > 0.\) For~\(\epsilon > 0\) we therefore get from Chebychev's inequality that
\[\mathbb{P}_{\delta}\left(|\beta_{tot}| > \epsilon \left(\frac{n}{A^2}+1\right) \mathbb{E}_{\delta}\tilde{T}_1^{\alpha}\right) \leq D_1\frac{\mathbb{E}_{\delta}(\beta_{tot}^m)}{n^m\epsilon^m} \leq \frac{D_2}{n^{m/2}\epsilon^m}\]
for some constants~\(D_1,D_2 > 0.\) Setting~\(\epsilon= \frac{1}{n^{1/16}}\) and using~\(\epsilon \left(\frac{n}{A^2}+1\right)  \leq \left(1+\frac{1}{n^{1/16}}\right) \frac{n}{A^2}\) for all~\(n\) large, we then get
\begin{equation}\label{til_est1}
\mathbb{P}_{\delta}\left(\sum_{i=1}^{\frac{n}{A^2}+1} \tilde{T}^{\alpha}_i \leq  \left(1+\frac{1}{n^{1/16}}\right)\frac{n}{A^2} \mathbb{E}\tilde{T}_1^{\alpha}\right) \geq 1-\frac{D_2}{n^{7m/16}}.
\end{equation}

From the upper bound for~\(S^{(\delta)}_{\alpha}\) in~(\ref{dom1}) and the fact that there is at least one node of~\({\cal P}_{\delta}\) in the unit square~\(S\) with probability~\(1-e^{-\delta n}\) we further get
\begin{equation}
\mathbb{P}_{\delta}\left(S^{(\delta)}_{\alpha} \leq  \left(1+\frac{1}{n^{1/16}}\right)\frac{n}{A^2} \mathbb{E}\tilde{T}_1^{\alpha}\right) \geq 1-\frac{D_2}{n^{7m/16}}- e^{-\delta n} \geq 1-\frac{2D_2}{n^{7m/16}}\nonumber
\end{equation}
for all~\(n\) large. Using the coupling relation~(\ref{mon_salpha}) we finally get~(\ref{til_est2}) proving Lemma~\ref{s_alpha_lemma}.~\(\qed\)

%PRKVMM ABV +eTC..FOR OUR BNFT.. +etC..

\setcounter{equation}{0}
\renewcommand\theequation{\thesection.\arabic{equation}}
\section{Variance upper bound for the MST}\label{mst_var_up}
In this section, we study variance upper bound estimates for MSTs with location dependent edge weights. Recalling the distribution parameters~\(\epsilon_1,\epsilon_2\) (see~(\ref{f_eq})) and the edge weight distribution parameters~\(c_1,c_2\) (see~(\ref{eq_met1}), we have the following result regarding the variance of~\(MST_n\) as defined in~(\ref{min_weight_tree}).
\begin{Theorem}\label{thm_mst_var_up} There is a constant~\(D_1 = D_1(\alpha,c_1,c_2,\epsilon_1,\epsilon_2) >  0\) such that the variance~\(var(MST_n) \leq D_1 n^{1-\alpha}\) for all~\(n\) large.
\end{Theorem}
As before, we have that the variance upper bound in the location dependent case is of the same order as that of the location independent case~(Kesten and Lee (1996)).

\underline{\emph{Remarks}}: We derive the variance estimate in Theorem~\ref{thm_mst_var_up} using the standard procedure (see Steele (1988), Kesten and Lee (1996)) of first obtaining one node difference estimates and then using the martingale difference method. Because of the large degree property described in Proposition~\ref{prop1}, we use a slightly different method in obtaining one node difference estimates in the location dependent case. We explain this in more detail in the paragraph following Lemma~\ref{f_xj_est} below.

\subsection*{\em One node difference estimates}
As in~(\ref{min_weight_tree}), let~\(MST_{n+1} = W({\cal T}_{n+1})\) be the weight of the MST~\({\cal T}_{n+1}\) formed by the nodes~\(\{X_k\}_{1 \leq k \leq n+1}\) and for~\(1 \leq j \leq n+1,\) let~\(MST_n(j) = W({\cal T}_n(j))\) be the weight of the MST~\({\cal T}_{n}(j)\) formed by the nodes~\(\{X_k\}_{1 \leq k \neq j \leq n+1}.\) We are interested in estimates for~\(|MST_{n+1} - MST_n(j)|,\) the change in the weight of the MST upon adding or removing a single node.

Consider the MST~\({\cal T}_{n}(j)\) with vertex set~\(\{X_i\}_{1 \leq i \neq j \leq n+1}\) and suppose~\(X_{i_0},\)\\\(1 \leq i_0 \neq j \leq n+1\) is the node closest to~\(X_j\) in terms of the Euclidean distance. The union~\({\cal T}_{n}(j) \cup \{(X_{j},X_{i_0})\}\) is a spanning tree containing all the nodes~\(\{X_k\}_{1 \leq k \leq n+1}\) and has weight~\(MST_{n}(j) + h^{\alpha}(X_{j},X_{i_0}).\) Therefore\\\(MST_{n+1} \leq MST_{n}(j) + h^{\alpha}(X_{j},X_{i_0})\) and using the bounds for the weight function~\(h(.)\) in~(\ref{eq_met1}), we have~\[h(X_{j},X_{i_0}) \leq c_2 d(X_{j},X_{i_0}) = c_2 d\left(X_j,\{X_k\}_{1 \leq k \neq j \leq n+1}\right),\] where~\(d(X_j,\{X_k\}_{1 \leq k \neq j \leq n+1})\) denotes the minimum distance between~\(X_j\) and the rest of the nodes. Thus
\begin{equation}\label{mst_up1}
MST_{n+1} \leq MST_{n}(j) + c_2^{\alpha}d^{\alpha}\left(X_j,\{X_k\}_{1 \leq k \neq j \leq n+1}\right).
\end{equation}

For getting an estimate in the reverse direction, we let~\(d_j\) be the degree of the node~\(X_j\) in the MST~\({\cal T}_{n+1}\) and let~\({\cal N}(X_j,{\cal T}_{n+1}) = \{v_1,\ldots,v_{d_j}\}\) be the set of neighbours of~\(X_j\) in~\({\cal T}_{n+1}.\) We remove the node~\(X_j\) and add the edges~\((v_i,v_{i+1}),1 \leq i \leq d_j-1\) to get a spanning tree containing all the nodes~\(\{X_k\}_{1 \leq k \neq j \leq n+1}.\) Again using~\(h(v_i,v_{i+1}) \leq c_2 d(v_i,v_{i+1})\) from~(\ref{eq_met1}), we get
\begin{equation}\label{mst_up2}
MST_n(j) \leq MST_{n+1} + \sum_{i=1}^{d_j-1}h^{\alpha}(v_i,v_{i+1}) \leq MST_{n+1} + c_2^{\alpha}\sum_{i=1}^{d_j-1}d^{\alpha}(v_i,v_{i+1})
\end{equation}
For~\(1 \leq i \leq d_j-1,\) we have by triangle inequality that~\(d(v_i,v_{i+1}) \leq d(X_j,v_i) + d(X_j,v_{i+1})\) and so using~\((a+b)^{\alpha} \leq 2^{\alpha}(a^{\alpha}+b^{\alpha})\) for all~\(a,b,\alpha > 0\) we get that~\(d^{\alpha}(v_i,v_{i+1}) \leq 2^{\alpha}\left(d^{\alpha}(X_j,v_i) + d^{\alpha}(X_j,v_{i+1})\right).\) Using this estimate in~(\ref{mst_up2}), we get
\begin{eqnarray}
MST_n(j) &\leq& MST_{n+1} + (2c_2)^{\alpha} \sum_{i=1}^{d_j}d^{\alpha}(X_j,v_{i}) \nonumber\\
&=& MST_{n+1} + (2c_2)^{\alpha} \sum_{v \in {\cal N}(X_j,{\cal T}_{n+1})}d^{\alpha}(X_j,v). \label{mst_up3}
\end{eqnarray}
Summarizing we get from~(\ref{mst_up1}) and~(\ref{mst_up3}) that
\begin{equation}\label{cruc_est}
|MST_{n+1} - MST_n(j)| \leq f_1(X_j)+f_2(X_j),
\end{equation}
where~\[f_1(X_j)  := c_2^{\alpha} d^{\alpha}\left(X_j,\{X_k\}_{1 \leq k \neq j \leq n+1}\right)\] and
\begin{equation}\label{f2_j_def}
f_2(X_j) := (2c_2)^{\alpha}\sum_{v \in {\cal N}(X_j,{\cal T}_{n+1})}d^{\alpha}(X_j,v).
\end{equation}

For future use, we prove the following Lemma.
\begin{Lemma}\label{f_xj_est}
There is a constant~\(D > 0\) such that for all~\(n\) large and any~\(1 \leq j \leq n+1\) we have
\begin{eqnarray}\label{f1_prop1}
\left(\mathbb{E}f_1(X_j)\right)^2 &\leq& \mathbb{E}f_1^2(X_j) \leq \frac{D}{n^{\alpha}},\nonumber\\
\sum_{j=1}^{n+1} \mathbb{E}f_1^2(X_j) &=& (n+1)\mathbb{E}f_1^2(X_1) \leq 2Dn^{1-\alpha},
\end{eqnarray}
\begin{eqnarray}\label{f2_prop1}
f_2^2(X_j) &\leq& D\sum_{v \in {\cal N}(X_j,{\cal T}_{n+1})}d^{2\alpha}(X_j,v),\nonumber\\
\sum_{j=1}^{n+1} \mathbb{E}f_2^2(X_j) &=& (n+1)\mathbb{E}f_2^2(X_1) \leq 2Dn^{1-\alpha}
\end{eqnarray}
and
\begin{equation}\label{one_diff_est_mst}
\mathbb{E}|MST_{n+1} - MST_n| \leq \mathbb{E}f_1(X_n) + \mathbb{E}f_2(X_n) \leq \left(\frac{D}{\sqrt{n}}\right)^{\alpha}.
\end{equation}
\end{Lemma}
The proof of~(\ref{one_diff_est_mst}) follows from~(\ref{f1_prop1}) and~(\ref{f2_prop1}). Estimate~(\ref{f1_prop1}) is standard and does not depend on the structure of the MST since it requires estimating the minimum distance between~\(X_j\) and the rest of the nodes. In case the edge weights are location independent, estimate~(\ref{f2_prop1}) is true (see Kesten and Lee (1996)) since the maximum degree of any node of MST is at most~\(6\) (see Aldous and Steele (1992) for example) and so  the first estimate in~(\ref{f2_prop1}) follows from the definition of~\(f_2\) in~(\ref{f2_j_def}) and the identity~\(\left(\sum_{i=1}^{6} a_i\right)^2 \leq 6 \sum_{i=1}^{6} a_i^2.\) The second estimate in~(\ref{f2_prop1}) then follows from the~\(\alpha-\)invariance property of MST (Kesten and Lee (1996)) which states that the MST remains the same whatever the value of the edge weight exponent~\(\alpha.\) This invariance property is also a consequence of the Kruskal's construction of the minimum spanning tree.

In our case where the edge weights are location dependent, we directly estimate the sum~\(f_2(X_j)\) by showing that the ratio of length of closely spaced edges having a common vertex is bounded above by a constant strictly less than one. Splitting the plane into sectors then allows us to estimate the sum of weighted length of edges within each sector and then use the~\(\alpha-\)invariance property as described above to obtain~(\ref{f2_prop1}) (see proof of~(\ref{f2_prop1}) below) .

For completeness we begin with the proof of~(\ref{f1_prop1}).\\
\emph{Proof of~(\ref{f1_prop1}) in Lemma~\ref{f_xj_est}}: Letting~\(d_{min}(X_j) = d\left(X_j,\{X_k\}_{1 \leq k \neq j \leq n+1}\right)\) we condition on~\(X_j~=~x\) and get from Fubini's theorem that
\begin{equation}\label{fub2}
\mathbb{E}f_1^2(X_j) = c_2^{2\alpha}\mathbb{E}d^{2\alpha}_{min}(X_j) = c_2^{2\alpha} \int \mathbb{E}d^{2\alpha}_{min}(x)f(x)dx,
\end{equation}
where
\begin{equation}\label{fub}
\mathbb{E}d^{2\alpha}_{min}(x) = 2\alpha \int y^{2\alpha-1}\mathbb{P}\left(d_{min}(x) > y\right) dy.
\end{equation}

The minimum distance from~\(x\) to~\(\{X_k\}_{1 \leq k \neq j \leq n+1}\) is at least~\(y\) if and only if~\(B(x,y) \cap S\) contains no node of~\(\{X_k\}_{1 \leq k \neq j \leq n+1},\) where~\(B(x,y)\) is the ball of radius~\(y\) centred at~\(x\) and we recall that~\(S\) is the unit square. The area of~\(B(x,y) \cap S\) is at least~\(\frac{\pi y^2}{4}\) no matter the location of~\(x\) and so using~\(f(x) \geq \epsilon_1\) from~(\ref{f_eq}), we get
\begin{equation}
\mathbb{P}(d_{min}(x) > y) = \left(1-\int_{B(x,y) \cap S} f(x)dx\right)^{n} \leq \left(1-\frac{\epsilon_1 \pi y^2}{4} \right)^{n} \leq e^{-\frac{n\epsilon_1 \pi y^2}{4}}. \nonumber
\end{equation}
For integer~\(l \geq 0\) we therefore have that~\(\mathbb{P}(\sqrt{n} \cdot d_{min}(x) > l) \leq e^{-C \cdot l^2}\) for some constant~\(C > 0\)
not depending on~\(X.\)
Thus
\begin{equation}\label{fubb3}
\mathbb{E}\left(\sqrt{n}\cdot d_{min}(x) \right)^{\alpha} \leq \sum_{l \geq 0}\mathbb{P}\left(\sqrt{n} \cdot d_{min}(x)\geq l^{\frac{1}{\alpha}}\right) \leq \sum_{l \geq 0} e^{-C \cdot l^{\frac{2}{\alpha}}},
\end{equation}
where the final summation in~(\ref{fubb3}) does not depend on~\(x.\) Substituting~(\ref{fubb3}) into~(\ref{fub2}) completes  the proof of~(\ref{f1_prop1}).~\(\qed\)

\emph{Proof of~(\ref{f2_prop1}) in Lemma~\ref{f_xj_est}}: To prove the first estimate in~(\ref{f2_prop1}), we let~\(K\) be a large integer satisfying
\[\left(2-2\cos\left(\frac{2\pi}{K}\right)\right)^{\frac{1}{2}}  < \frac{c_1}{c_2} < 1,\] where~\(c_1\) and~\(c_2\)
are the bounds for the weight function~\(h\) as in~(\ref{eq_met1}). Defining~\[g(x) := \left(1+x^2 - 2\cdot x\cdot \cos\left(\frac{2\pi}{K}\right)\right)^{\frac{1}{2}}, 0 \leq x \leq 1\] we have that~\(g(0) = 1\)
and since~\(g\) is continuous in~\([0,1],\) there exists a positive number~\(r_0 < 1\) such that~\(g(x) \geq \frac{c_1}{c_2}\)
for~\(0 < x \leq r_0\) and~\(g(x) < \frac{c_1}{c_2}\) for\\\(r_0 < x \leq 1.\)

%Let~\(\{v_1,\ldots,v_{d_j}\}\) be the neighbours of~\(X_j\) in the spanning tree~\({\cal T}_{n+1}\) encountered in the clockwise order
%and for convenience set~\(v_0 := v_{d_j}\) and~\(v_{d_j+1} := v_1.\) For~\(1 \leq l \leq d_j\) let~\(\theta_{l}\) denote the angle between the edges~\((X_j,v_l)\) and~\((X_j,v_{l+1}).\)

Draw~\(K\) rays starting from~\(X_j\) equally spaced at an angle of~\(\frac{2\pi}{K}\) apart; i.e.,
let~\(l_1,\ldots,l_K\) be~\(K\) rays starting from~\(X_j\) encountered in the clockwise order such that the angle between~\(l_i\) and~\(l_{i+1}\) is~\(\frac{2\pi}{K}.\) Let~\(v_{i_1},\ldots,v_{i_w}\) be the neighbours of~\(X_j\) present in the~\(i^{th}\) sector formed by the rays~\(l_i\) and~\(l_{i+1}\) and suppose that~\(d(X_j,v_{i_1}) > d(X_j,v_{i_2}) > \ldots > d(X_j,v_{i_w}).\) We first show that the edge length ratio~\(r := \frac{d(X_j,v_{i_{k+1}})}{d(X_j,v_{i_{k}})} \leq r_0\) for~\(1 \leq k \leq w-1\) where~\(r_0\) is as in the previous paragraph.

If~\(\theta\) denotes the angle between the edges~\((X_j,v_{i_k})\) and~\((X_j,v_{i_{k+1}}),\) then\\\(\theta < \frac{2\pi}{K}\) and so
\begin{eqnarray}
d^2(v_{i_k},v_{i_{k+1}}) &=& d^2(X_j,v_{i_k}) + d^2(X_j,v_{i_{k+1}}) - 2d(X_j,v_{i_k})\cdot d(X_j,v_{i_{k+1}})\cdot \cos(\theta) \nonumber\\
&=& d^2(X_j,v_{i_k})(1+r^2-2r\cos(\theta)) \nonumber\\
&\leq& d^2(X_j,v_{i_k})\left(1+r^2-2r\cos\left(\frac{2\pi}{K}\right)\right). \nonumber
\end{eqnarray}
Thus~\(d(v_{i_k},v_{i_{k+1}}) \leq d(X_j,v_{i_k}) g(r)\) and if~\(r_0 < r \leq 1,\) then using~(\ref{eq_met1}) we have~\[ h(v_{i_k},v_{i_{k+1}}) \leq c_2 d(v_{i_k},v_{i_{k+1}}) \leq c_2 d(X_j,v_{i_k}) g(r) < c_1d(X_j,v_{i_k})\leq h(X_j,v_{i_k}),\]
by our choice of~\(r_0\) in the first paragraph. This is a contradiction because removing the edge~\((X_j,v_{i_k})\)
from the tree~\({\cal T}_{n+1}\) and adding the edge~\((v_{i_k},v_{i_{k+1}}),\) we get a spanning tree with weight less than~\(MST_{n+1}.\)

Using the ratio property iteratively, we get~\(d(X_j,v_{i_k}) \leq r_0^{k-1} d(X_j,v_{i_1})\)
and so~\(\sum_{k=1}^{w} d^{\alpha}(X_j,v_{i_k}) \leq C d^{\alpha}(X_j,v_{i_1}),\) where~\(C = \sum_{k \geq 1} r_0^{\alpha(k-1)}.\)
This estimate holds for each of the~\(K\) sectors and so
\[\sum_{v \in {\cal N}_{n+1}(X_j,{\cal T}_{n+1})}d^{\alpha}(X_j,v) \leq \sum_{i=1}^{K} \sum_{v} d^{\alpha}(X_j,v) \leq C \sum_{i=1}^{K}d^{\alpha}(X_j,u_i),\]
where the middle summation is over all nodes~\(v\) present in sector~\(i\) and~\(u_i\) in the final summation is the node farthest in Euclidean distance from~\(X_j\) in the~\(i^{th}\) sector. Therefore using~\((\sum_{i=1}^{K}b_i)^2 \leq K \sum_{i=1}^{K} b_i^2\) for positive numbers~\(\{b_i\}\) we get
\begin{eqnarray}
\left(\sum_{v \in {\cal N}_{n+1}(X_j,{\cal T}_{n+1})}d^{\alpha}(X_j,v)\right)^{2} &\leq& C^2 \left(\sum_{i=1}^{K}d^{\alpha}(X_j,u_i)\right)^2 \nonumber\\
&\leq& C^2 K \sum_{i=1}^{K}d^{2\alpha}(X_j,u_i) \nonumber\\
&\leq& C^2 K \sum_{v \in {\cal N}_{n+1}(X_j,{\cal T}_{n+1})} d^{2\alpha}(X_j,v), \nonumber
\end{eqnarray}
proving the first estimate in~(\ref{f2_prop1}).

%WRT HERE THAT THE MST REMAINS THE SAME WHATEVER THE VALUE OF ALPHA...

Consequently~
\begin{eqnarray}
\sum_{j=1}^{n+1} \mathbb{E}f_2^2(X_j) &\leq& C\mathbb{E}\sum_{j=1}^{n+1} \sum_{v \in {\cal N}(X_j,{\cal T}_{n+1})} d^{2\alpha}(X_j,v) \nonumber\\
&\leq& \frac{C}{c_1^{2\alpha}}\mathbb{E}\sum_{j=1}^{n+1} \sum_{v \in {\cal N}(X_j,{\cal T}_{n+1})} h^{2\alpha}(X_j,v), \label{pun_sum}
\end{eqnarray}
again using~(\ref{eq_met1}). The double summation in the final term is simply twice the weight of the MST with edge weight exponent~\(2\alpha.\) This is because, given the location of the nodes~\(\{X_i\}_{1 \leq i \leq n+1},\) the edge weights~\(h(X_i,X_j)\) are fixed and so the Kruskal's algorithm gives the same MST irrespective of the value of the edge weight exponent~\(\alpha\) (Kesten and Lee (1996), Yukich (2000)); i.e., if we denote~\( {\cal T}_{n+1}(\beta)\) to be the MST for edge weight exponent~\(\beta > 0\) as in~(\ref{min_weight_tree}) so that \[W({\cal T}_{n+1}(\beta)) = \sum_{e \in {\cal T}_{n+1}(\beta)} h^{\beta}(e) = \min_{\cal T} W({\cal T}) = \min_{\cal T} \sum_{f \in {\cal T}}h^{\beta}(f),\] where the minimum is taken over all spanning trees~\({\cal T}\) containing all the~\(n~+~1\) nodes~\(\{X_i\}_{1 \leq i \leq n+1},\) then~\({\cal T}_{n+1}(\beta) = {\cal T}_{n+1}(1)\) for any~\(\beta > 0.\) Therefore using~(\ref{pun_sum}) and the expectation upper bound~(\ref{exp_mst_bound}) with~\(2\alpha\) instead of~\(\alpha,\) we also get the second estimate in~(\ref{f2_prop1}).~\(\qed\)

\subsection*{\em Proof of Theorem~\ref{thm_mst_var_up}}
We use one node difference estimate~(\ref{cruc_est}) together with the martingale difference method to obtain a bound for the variance. For~\(1 \leq j \leq n~+~1,\) let~\({\cal F}_j = \sigma\left(\{X_k\}_{1 \leq k \leq j}\right)\) denote the sigma field
generated by the node positions~\(\{X_k\}_{1 \leq k \leq j}.\) Defining the martingale difference
\begin{equation}\label{h_def}
H_j = \mathbb{E}(MST_{n+1} | {\cal F}_j) - \mathbb{E}(MST_{n+1} | {\cal F}_{j-1}),
\end{equation}
we then have that~\(MST_{n+1} -\mathbb{E}MST_{n+1} = \sum_{j=1}^{n+1} H_j\) and so by the martingale property
\begin{equation} \label{var_exp}
var(MST_{n+1})  = \mathbb{E}\left(\sum_{j=1}^{n+1} H_j\right)^2 = \sum_{j=1}^{n+1} \mathbb{E}H_j^2.
\end{equation}

To evaluate~\(\mathbb{E}H_j^2\) we rewrite the martingale difference~\(H_j\) in a more convenient form.
Letting~\(X'_j\) be independent copy of~\(X_j\) which is also independent of~\(\{X_k\}_{1 \leq k \neq j \leq n+1}\) we rewrite~
\begin{equation}\label{h_def2}
H_j = \mathbb{E}(MST_{n+1}(X_j) - MST_{n+1}(X'_j) | {\cal F}_j),
\end{equation}
where~\(MST_{n+1}(X_j)\) is the weight of the MST formed by the nodes~\(\{X_i\}_{1 \leq i \leq n+1}\) and~\(MST_{n+1}(X'_j)\)
is the weight of the MST formed by the nodes\\\(\{X_i\}_{1 \leq i \neq j \leq n+1}~\cup~\{X'_j\}.\)
Using the triangle inequality and the one node difference estimate~(\ref{cruc_est}), we have
that~\(|MST_{n+1}(X_j)-MST_{n+1}(X'_j)|\) is bounded above as
\begin{eqnarray}
&&|MST_{n+1}(X_j)-MST_{n}(j)| + |MST_{n+1}(X'_j)-MST_{n}(j)| \nonumber\\
&&\;\;\;\leq\;\;f_1(X_j) + f_2(X_j) + f_1(X'_j) + f_2(X'_j), \nonumber
\end{eqnarray}
where~\(f_1\) and~\(f_2\) are as in~(\ref{cruc_est}) and we recall that~\(MST_n(j)\) is the weight of the MST
formed by the nodes~\(\{X_k\}_{1 \leq k \neq j \leq n+1}.\)
Thus
\begin{eqnarray}
|H_j| &\leq& \mathbb{E}(|MST_{n+1}(X_j)-MST_{n+1}(X'_j)| | {\cal F}_j) \nonumber\\
&\leq& \mathbb{E}(f_1(X_j)|{\cal F}_j) + \mathbb{E}(f_2(X_j) |{\cal F}_j) + \mathbb{E}(f_1(X'_j)|{\cal F}_{j}) + \mathbb{E}(f_2(X'_j)|{\cal F}_{j}) \nonumber\\
&=& \mathbb{E}(f_1|{\cal F}_j) + \mathbb{E}(f_2 |{\cal F}_j) + \mathbb{E}(f_1|{\cal F}_{j-1}) + \mathbb{E}(f_2|{\cal F}_{j-1}). \nonumber
\end{eqnarray}
Using~\((a_1+a_2+a_3+a_4)^2 \leq 4(a_1^2 + a_2^2 + a_3^2 + a_4^2),\) we then get
\begin{eqnarray}
H_j^2 &\leq& 4\left(\left(\mathbb{E}(f_1|{\cal F}_j)\right)^2 + \left(\mathbb{E}(f_2 |{\cal F}_j)\right)^2 + \left(\mathbb{E}(f_1|{\cal F}_{j-1})\right)^2 + \left(\mathbb{E}(f_2|{\cal F}_{j-1})\right)^2\right)  \nonumber\\
&\leq& 4\left(\mathbb{E}(f^2_1|{\cal F}_j) + \mathbb{E}(f^2_2 |{\cal F}_j)+ \mathbb{E}(f^2_1|{\cal F}_{j-1}) + \mathbb{E}(f^2_2|{\cal F}_{j-1})\right) \nonumber
\end{eqnarray}
since~\((\mathbb{E}(X|{\cal F}))^2 \leq \mathbb{E}(X^2|{\cal F}).\) Thus~\(\mathbb{E}H_j^2 \leq 8\left(\mathbb{E}f_1^2(X_j) + \mathbb{E}f_2^2(X_j)\right)\) and plugging this in~(\ref{var_exp}), we have
\[var(MST_{n+1}) = \sum_{j=1}^{n+1} \mathbb{E}H_j^2 \leq  8\left(\sum_{j=1}^{n+1}\mathbb{E}f_1^2(X_j) + \sum_{j=1}^{n+1} \mathbb{E}f_2^2(X_j)\right),\]
which in turn is at most a constant multiple of~\(n^{1-\alpha}\) using the estimates in~(\ref{f1_prop1}) and~(\ref{f2_prop1}).~\(\qed\)

%SEE MORE +etC..

%WRT SMALL PROF OF COR +etC...

%In the end, we convert the estimates to the Binomial process.

%As in the case of MST, the variance upper bounds are obtained via the martingale difference method that estimates the change in TSP lengths after adding or removing a single node. For the variance lower bound of MST, we study the probability of occurrence of predetermined ``nice" configurations and thereby estimate the change in the MST length upon shifting a single node (see Section~\ref{pf_mst_var_low}).

%SEE ABV +eTC..

%do we write uniform separately??

%CHNG HERE.....The paper is organized as follows. WRT HERE...

%CHK REF ETC.. AND ALL +eTC...

%WRT SMALL PROF OF COR +etC...

%In the end, we convert the estimates to the Binomial process.

%Arguing iteratively as in the proof of~\((h2),\) we get
%\begin{eqnarray}
%\left|MST_{k} - MST_{n^2}\right| &\leq& C_1r_{n^2} (\log{n^2})(k-n^2)\ind\left(Y_{tot}(n^2)\right) \nonumber\\
%&&\;\;\;\;\;\;+\;\;\;(k-n^2) n^2\sqrt{2} \ind\left(Y_{tot}^c(n^2)\right)\label{d_est_mst}
%\end{eqnarray}
%GOOTHALS HERE...

\setcounter{equation}{0}
\renewcommand\theequation{\thesection.\arabic{equation}}
\section{Variance lower bound for MST}\label{mst_var_low}
Regarding the variance lower bound for~\(MST_n\) as defined in~(\ref{min_weight_tree}), we have the following result.
\begin{Theorem}\label{thm_mst_var_low}
Suppose there is a square~\(S_0 \subseteq S\) with constant side length~\(s_0\) such that the weight function~\(h(u,v) = d(u,v)\) if either~\(u\) or~\(v\) is in~\(S_0.\) There is a constant~\(D_2 = D_2(\alpha,c_1,c_2,\epsilon_1,\epsilon_2) >0\) such that~\(var(MST_n) \geq D_2 n^{1-\alpha}\) for all~\(n\) large.
\end{Theorem}
Combining with Theorem~\ref{thm_mst_var_up}, we then find that~\(var(MST_n)\) is of the \emph{order} of~\(n^{1-\alpha}\) even in the location dependent case. As before, the behaviour is analogous as in the location independent case and in this aspect we refer to Kesten and Lee (1996) who use martingale methods to study central limit theorems for~\(MST_n,\) appropriately scaled and centred. The technical assumption of weight function being equal to the Euclidean length in a small subsquare within the unit square allows us to bound probabilities of predetermined nice configurations via martingale difference estimates (see proof of Theorem~\ref{thm_mst_var_low} below).

We begin with some preliminary definitions and computations. Suppose that~\(s_0\) is the side length of the square~\(S_0.\) As in the proof of the deviation estimates in Section~\ref{mst_dev}, we divide~\(S_0\) into small~\(\frac{A}{\sqrt{n}} \times \frac{A}{\sqrt{n}}\) squares~\(\{R_i\}_{1 \leq i \leq \frac{s^2_0n}{A^2}}\) where~\(A  = A(n)  \in \left[1, 1+\frac{1}{\log{n}}\right]\) is chosen such that~\(\frac{s_0\sqrt{n}}{A}\) is an integer. This is possible by the argument preceding~(\ref{poss}).

\begin{figure}
\centering
\begin{subfigure}{0.5\textwidth}
\centering
\includegraphics[width=2.5in, trim= 50 200 50 200, clip=true]{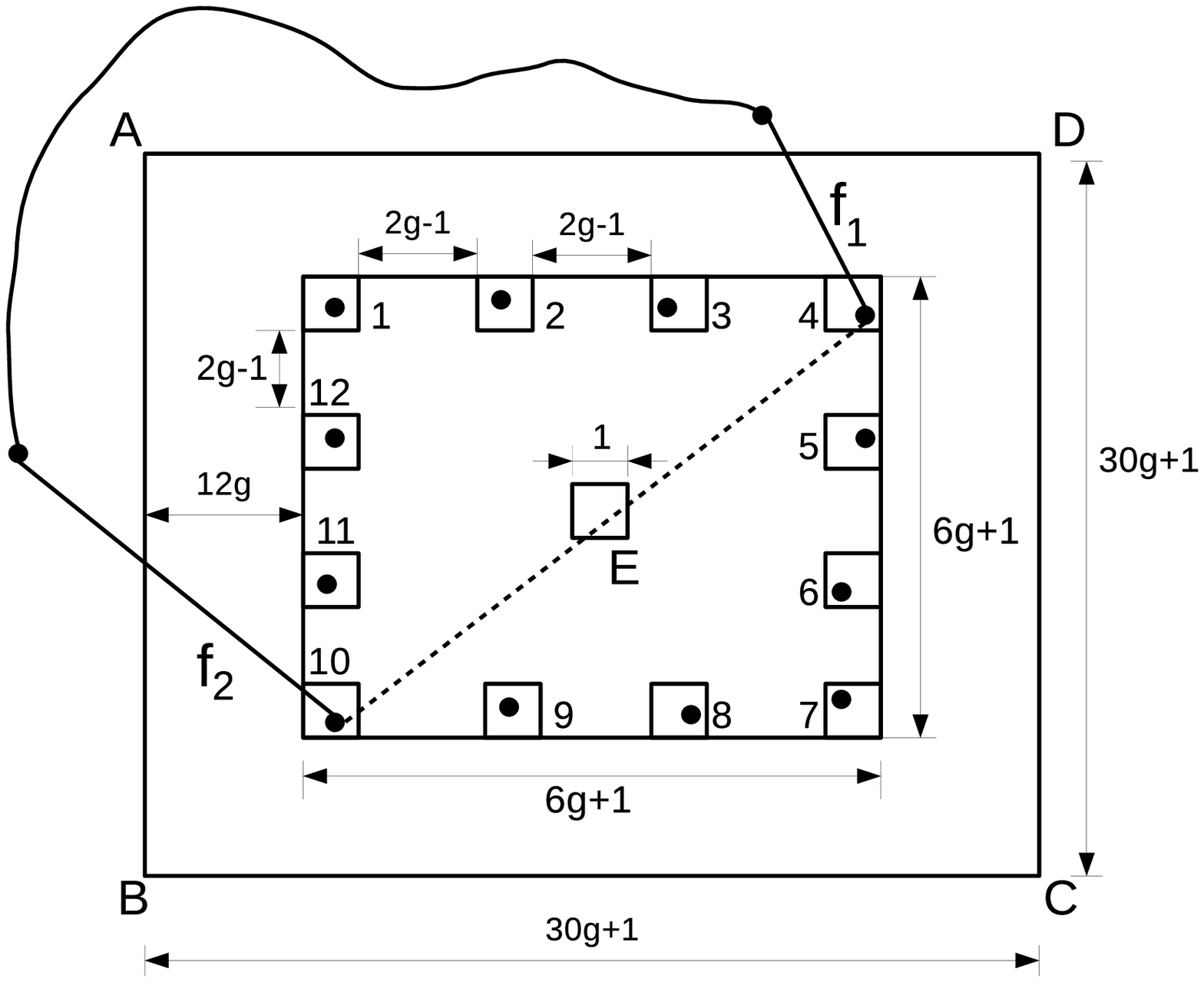}
  \caption{} %\caption{Tiling~\(S\) into~\(r_n \times r_n\) squares with an inter-square distance of~\(s_n.\)}%\nonumber%{t1}
\end{subfigure}% seems this is important....
\begin{subfigure}{0.5\textwidth}
\centering
   \includegraphics[width=2.5in, trim= 50 200 50 200, clip=true]{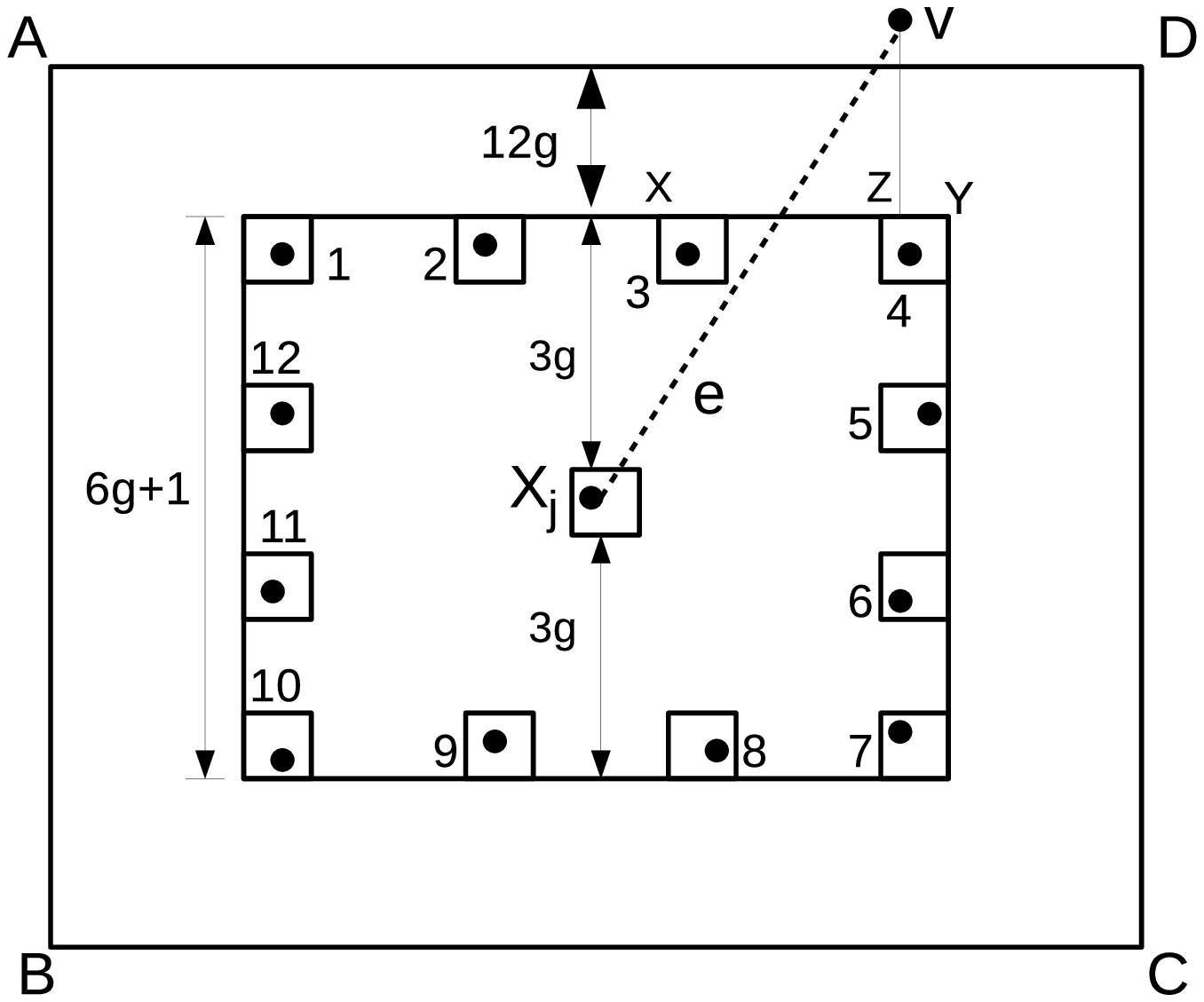}
   \caption{} %\caption{The induced subgraph~\(G({\cal D}_n(N))\) corresponding to the centres of the grey squares in~\((a).\)} %\nonumber%{t2}
  \end{subfigure}

\caption{\((a)\) The square~\(E\) is a ~\(g-\)good square. \((b)\) The edge~\(e\) has longer Euclidean length than the edge~\((v,v_3).\) All dimensions are to be multiplied by~\(\frac{A}{\sqrt{n}}\) to get the actual length.}
\label{good_square}
\end{figure}

For a square~\(R_i\) we let~\({\cal N}_1(R_i)\) be the squares in~\(\{R_k\}\) sharing a corner with~\(R_i\) and for~\(l \geq 2\) we let~\({\cal N}_{l}(R_j)\) be the set of all squares in~\(\{R_k\}\) sharing a corner with~\({\cal N}_{l-1}(R_j).\) Thus~\({\cal N}_l(R_i)\) contains~\((2l+1)^2\)
of~\(\{R_k\}\) including the square~\(R_i.\) For an integer~\(g \geq 5,\) we depict~\({\cal N}_{3g}(R_i)\) and~\({\cal N}_{15g}(R_i)\) in Figure~\ref{good_square}~\((a)\) as the~\((6g+1) \times (6g+1) \) square and~\((30g+1) \times (30g+1)\) square, respectively. In all the figures in this subsection, the dimensions are to be multiplied by the scaling factor~\(\frac{A}{\sqrt{n}}.\)
The central square labelled~\(E\) is the~\(\frac{A}{\sqrt{n}} \times \frac{A}{\sqrt{n}}\) square~\(R_i\) and the twelve~\(\frac{A}{\sqrt{n}} \times \frac{A}{\sqrt{n}}\) squares numbered~\(1,2,\ldots,12\)  are all spaced~\((2g-1)\frac{A}{\sqrt{n}}\) apart.

\begin{Definition}\label{def_gj}
For~\(1 \leq j \leq n+1\) we say that the square~\(R_i\) is a~\((g,j)-\)\emph{good square} if each of the~\(\frac{A}{\sqrt{n}} \times \frac{A}{\sqrt{n}}\) squares numbered one to twelve in Figure~\ref{good_square}~\((a)\) contain exactly one node of~\(\{X_k\}_{1 \leq k \neq j \leq n+1}\) and the rest of the big square~\(ABCD\) (containing all the~\(\frac{A}{\sqrt{n}} \times \frac{A}{\sqrt{n}}\) squares of~\({\cal N}_{15}(R_i)\)) contains no node of~\(\{X_k\}_{1 \leq k \neq j \leq n+1}.\)
\end{Definition}

The advantage of defining the~\((g,j)-\)good square as above is that we can determine exactly the change in the MST length when the node~\(X_j\) is ``added" to the~\((g,j)-\)good square~\(R_{i}.\) Indeed, suppose that\(R_i\) is a~\((g,j)-\)good square and the node~\(X_j \in R_i.\) For~\(1 \leq l \leq 12,\) let~\(v_l\) be the node present in the~\(\frac{A}{\sqrt{n}} \times \frac{A}{\sqrt{n}}\) square numbered~\(l\) in Figure~\ref{good_square} and let~\(v_{min} \in \{v_l\}_{1 \leq l \leq 12}\) be the node closest to~\(X_j\) in terms of Euclidean length. As before we denote~\({\cal T}_{n}(j)\) and~\({\cal T}_{n+1}\) to be the MSTs formed by the nodes~\(\{X_i\}_{1 \leq i \neq j \leq n+1}\) and~\(\{X_i\}_{1 \leq i \leq n+1},\) respectively, with edge weight function~\(h(.,.)\) and edge weight exponent~\(\alpha > 0.\) Letting~\(MST_n(j)\) denote the weight of~\({\cal T}_n(j),\) we have the following result.
\begin{Lemma}\label{var_low_lem1} If the event~\(\{R_i \text{ is }(g,j)-\text{good}\} \cap \{X_j \in R_i\} \) occurs then
\begin{equation}\label{add_one}
{\cal T}_{n+1} = {\cal T}_{n}(j) \cup \{(X_j,v_{min})\} =: {\cal T}_{new}.
\end{equation}
and
\begin{equation}\label{add_one_est1}
(3g-1)^{\alpha}\cdot \left(\frac{A}{\sqrt{n}}\right)^{\alpha} \leq MST_{n+1} - MST_n(j) \leq (5g-1)^{\alpha} \cdot \left(\frac{A}{\sqrt{n}}\right)^{\alpha}.
\end{equation}
\end{Lemma}
In words, the new MST~\({\cal T}_{n+1}\) is simply obtained by adding the edge~\((X_j,v_{min})\) to
the old MST~\({\cal T}_{n}(j).\) Estimate~(\ref{add_one_est1}) obtains explicit bounds on the difference of MST weight upon adding or removing the node~\(X_j.\) This allows us to use martingale difference estimates and obtain a lower bound on the variance in the proof of Theorem~\ref{thm_mst_var_low} below.

Before we do so, we state the following companion Lemma that shows that the event mentioned in the statement of Lemma~\ref{var_low_lem1} occurs with positive probability. For~\(1 \leq j \leq n+1\) let~\(X'_j\) be an independent copy of~\(X_j\) that is also independent of~\(\{X_k\}_{1 \leq k \neq j \leq n+1}\) and let~\(F_{j}(g)\) and~\(F'_j(g)\) respectively denote the events that the nodes~\(X_j\) and~\(X'_j\) belong to a~\((g,j)-\)good square in~\(\{R_k\}\)
and let~\(F^{tot}_j(g) = F_j(2g) \cap F'_j(g).\)
\begin{Lemma}\label{var_egj_lem} For every~\(g \geq 5,\) there is a constant~\(\theta = \theta(g,s_0,\epsilon_1,\epsilon_2) >0 \) such that \[\min_{1 \leq j \leq n+1} \mathbb{P}(F^{tot}_j(g)) \geq \theta.\]
\end{Lemma}
Assuming Lemmas~\ref{var_low_lem1} and~\ref{var_egj_lem}, we prove Theorem~\ref{thm_mst_var_low} below. We then prove Lemmas~\ref{var_low_lem1} and~\ref{var_egj_lem} separately.

\emph{Proof of Theorem~\ref{thm_mst_var_low} (assuming Lemmas~\ref{var_low_lem1} and~\ref{var_egj_lem})}: Letting~\(F^{tot}_j(g)\) be the event defined prior to Lemma~\ref{var_egj_lem}, we get from Lemma~\ref{var_low_lem1} that if~\(F^{tot}_j(g)\) occurs, then
\begin{eqnarray}
MST_{n+1}(X_j) &\geq& MST_{n}(j) + \left(\frac{A}{\sqrt{n}}\right)^{\alpha}(6g-1)^{\alpha} \nonumber\\
&\geq&  MST_{n+1}(X'_j) + \left(\frac{A}{\sqrt{n}}\right)^{\alpha}((6g-1)^{\alpha} -(5g-1)^{\alpha}) \nonumber
\end{eqnarray}
Letting~\(\Delta_{\alpha} := (6g-1)^{\alpha} -(5g-1)^{\alpha} > 0\)
and recalling the martingale difference~\(H_j\) in~(\ref{h_def2}), we get
\[\mathbb{E}|H_j| \geq \mathbb{E}|H_j|\ind(F^{tot}_j(g)) = \mathbb{E}H_j \ind(F^{tot}_j(g)) \geq \Delta_{\alpha}n^{-\frac{\alpha}{2}} \mathbb{P}(F^{tot}_j(g)) \geq c\Delta_{\alpha}n^{-\frac{\alpha}{2}},\] where~\(c > 0\) is the constant in Lemma~\ref{var_egj_lem}. Consequently from~(\ref{var_exp}), we get \[var(MST_{n+1}) = \sum_{j=1}^{n+1} \mathbb{E}H_j^2 \geq \sum_{j=1}^{n+1} \left(\mathbb{E}|H_j|\right)^2 \geq c^2\Delta_{\alpha}^2 n^{1-\alpha},\] proving the desired lower bound.~\(\qed\)

In the rest of this section, we prove Lemmas~\ref{var_low_lem1} and~\ref{var_egj_lem} starting with the former.

\subsection*{\em Proof of Lemma~\ref{var_low_lem1}}
To prove Lemma~\ref{var_low_lem1}, we denote the edges in~\(\{(v_i,v_{i+1})\}_{1 \leq i \leq 11} \cup \{(v_{12},v_1)\}\) to be \emph{short edges} and  collect the following properties, assuming that the event defined in Lemma~\ref{var_low_lem1} holds.\\
\((p1)\) The length of any short edge is at most~\((2g+5)\cdot\frac{A}{\sqrt{n}}.\) The subgraph~\({\cal T}_{loc}\) of~\({\cal T}_{n}(j)\) induced by the vertices~\(\{v_i\}_{1 \leq i \leq 12}\) is a tree consisting of exactly eleven short edges.\\
\((p2)\) Let~\(v\) be any point outside the square~\(ABCD\) as in Figure~\ref{good_square}\((b)\) so that the perpendicular~\(vZ\) from~\(v\) to the line containing~\(XY\) crosses the line segment~\(XY.\) The Euclidean distance~\(d(X_j,v) \geq 15g\cdot \frac{A}{\sqrt{n}}\) and~\(d(X_j,v) > d(v,v_3),\) strictly.\\

\emph{Proof of~\((p1)\)}: To prove that~\({\cal T}_{loc} \subset {\cal T}_n(j)\) is a tree we assume otherwise and suppose for example that nodes~\(v_{4}\) and~\(v_{10}\) are joined by a path~\(P_{4,10}\) in~\({\cal T}_{old}\) containing at least one vertex not in~\(\{v_i\}_{1 \leq i \leq 12}\) as in Figure~\ref{good_square}\((a).\) This means that~\(P_{4,10}\) contains at least two ``long" edges~\(f_1,f_2\) containing endvertices outside~\({\cal N}_{15g}(R_{u_1}).\)

The short edges and the edges~\(\{f_1,f_2\}\) all have an endvertex in the square~\(S_0\) and so the weights of these edges are simply their Euclidean lengths raised to the power~\(\alpha.\) The length of the edge~\(f_1\) is at least~\(12g \cdot \frac{A}{\sqrt{n}},\)  the width of the annulus between the big squares (see Figure~\ref{good_square}\((a)\)) but the length of the edge~\((v_4,v_{10})\) is at most~\((6g+1)\sqrt{2}\cdot\frac{A}{\sqrt{n}}.\) Since~\(g \geq 2\) we have~\((6g+1)\sqrt{2} < 10g\) and so removing the edge~\(f_1\) and adding the edge~\((v_4,v_{10}),\) we get an MST formed by the nodes~\(\{X_i\}_{1 \leq i \neq j \leq n+1}\) with weight strictly that~\({\cal T}_n(j),\) a contradiction. This proves that~\({\cal T}_{loc}\) is a tree and the edge set of~\({\cal T}_{loc}\) must contain only short edges, since any other edge with both endvertices in~\(\{v_i\}_{1 \leq i \leq 12}\) has length strictly larger than longest short edge.~\(\qed\)

\emph{Proof of~\((p2)\)}: From Figure~\ref{good_square}\((b),\) we have~\(d(X_j,v) \geq 15g\cdot \frac{A}{\sqrt{n}}.\) We show below that the angle~\(\theta\) between the edges~\((X_j,v)\) and~\((v,v_3)\) is less than~\(60\) degrees and so the Euclidean length~\(d(X_j,v) < d(v,v_3),\) strictly. To estimate~\(\theta\) we let~\(d(v,Z)\) be the length of the perpendicular segment~\(vZ\) so that the area of the triangle formed by the vertices~\(v,X\) and~\(Y\) is~\[\frac{1}{2}\cdot d(v,Z) \cdot d(X,Y) = \frac{1}{2}d(v,X)\cdot d(v,Y) \cdot \sin(\theta_0),\] where~\(\theta_0 > \theta\)
is the angle between the edges~\((v,X)\) and~\((v,Y).\) Thus~\(\sin(\theta_0) = \frac{d(v,Z)\cdot d(X,Y)}{d(v,X)\cdot d(v,Y)}\) and using~\(\min(d(v,X),d(v,Y)) \geq d(v,Z)\) we further get~\(\sin(\theta_0) \leq \frac{d(X,Y)}{d(v,Z)}.\) But~\(d(v,Z) \geq 12g\cdot\frac{A}{\sqrt{n}}\) and from Figure~\ref{good_square}\((a),\) we have~\(d(X,Y) \leq (2g+1)\cdot\frac{A}{\sqrt{n}}\) and so~\(\sin(\theta_0) \leq \frac{2g+1}{12g} \leq \frac{1}{4},\) since~\(g \geq 1.\)~\(\qed\)

\emph{Proof of~(\ref{add_one}) in Lemma~\ref{var_low_lem1}}: To prove~(\ref{add_one}), we first identify for every edge~\(e  = (u,v) \notin {\cal T}_{new},\) the unique path~\(P(e) \subseteq {\cal T}_{new}\) with endvertices~\(u,v.\) By Lemma~\(1\) of Kesten and Lee (1996), we then have that~\({\cal T}_{new}\) is the MST of the nodes~\(\{X_i\}_{1 \leq i \leq n+1}\) if and only if for every edge~\(e \notin {\cal T}_{new}\) the following holds:
\begin{equation}\label{mst_prop}
\text{for every edge~\(f \in P(e),\) the weight~\(h(f) < h(e)\)}
\end{equation}
where we recall that~\(h(e) = h(x,y)\) is the weight of the edge~\(e\) with endvertices~\(x\) and~\(y.\) The condition~(\ref{mst_prop}) does not depend on the value of the edge weight exponent~\(\alpha\) and as mentioned before, this obtains that the MST is the same irrespective of the value of~\(\alpha.\)

We now prove that~(\ref{mst_prop}) holds for each edge~\(e \notin {\cal T}_{new}.\) For notational convenience, we refer to~\({\cal T}_n(j)\) as~\({\cal T}_{old}.\) Suppose first that~\(e  = (u,v) \notin {\cal T}_{new}\) does not contain~\(X_j\) as an endvertex. There is a unique path~\(P(u,v) \subseteq {\cal T}_{new}\) with endvertices~\(u\) and~\(v\) and since~\(X_j\) is a leaf of~\({\cal T}_{new},\) no edge of~\(P(u,v)\) contains~\(X_j\) as an endvertex. Thus~\(P(u,v) \subseteq {\cal T}_{old}\) and so applying~(\ref{mst_prop}) to the MST~\({\cal T}_{old},\) we get that the weight of every edge in~\(P(u,v)\) is no more than the weight of~\(e,\) proving~(\ref{mst_prop}) for~\(e \in {\cal T}_{new}.\)

Suppose now that~\(e = (v,X_j)\) for some vertex~\(v\) as in Figure~\ref{good_square}\((b).\) The MST~\({\cal T}_{old}\) does not contain both the edges~\((v,v_3)\) and~\((v,v_4)\) since~\({\cal T}_{old}\) has eleven of the twelve short edges (property~\((p1)\)) and so if both~\((v,v_3)\) and~\((v,v_4)\) were to belong to~\({\cal T}_{old},\) this would create a cycle. Suppose~\((v,v_3) \notin {\cal T}_{old}\) so that there is a unique path~\(P(v,v_3) \subseteq {\cal T}_{old} \subset {\cal T}_{new}\) satisfying~(\ref{mst_prop}). Also let~\(P(v_3,v_{min}) \subseteq {\cal T}_{old} \subset {\cal T}_{new}\) be the unique path formed by the short edges with endvertices~\((v_3,v_{min})\) so that~\(P(v,v_3) \cup P(v_3,v_{min}) \cup \{(v_{min},X_j)\} \subseteq {\cal T}_{new}\) contains the unique path~\(P(v,X_j)\) in~\({\cal T}_{new}\) with endvertices~\(v\) and~\(X_j.\)

Applying~(\ref{mst_prop}) to the path~\(P(v,v_3) \subseteq {\cal T}_{old}\) we have that every edge in~\(P(v,v_3)\) has weight less than
\begin{equation}\label{cru_ew1}
h(v,v_3)  = d(v,v_3) < d(v,X_j) = h(v,X_j)
\end{equation}
where the first and the last equalities in~(\ref{cru_ew1}) follow from the fact that~\(v_3,X_j\) belong to the square~\(S_0.\)  The middle inequality in~(\ref{cru_ew1}) is true by property~\((p2).\) Also, every short edge has Euclidean length at most~\((2g+5)\cdot \frac{A}{\sqrt{n}}\) (property~\((p1)\)) and~\(d(v,X_j) \geq  15g\cdot \frac{A}{\sqrt{n}}\) (property~\((p2)\)) and so the weight of every edge in~\(P(v_3,v_{min})\) is also less than the weight of~\((X_j,v).\) Finally, the length of~\((X_j,v_{min})\) is at most the length of the diagonal of the inner big square in Figure~\ref{good_square}~\((a)\) and so~\(h(X_j,v_{min}) = d(X_{j},v_{min}) < (6g+1)\sqrt{2} \cdot \frac{A}{\sqrt{n}} < 15g \cdot \frac{A}{\sqrt{n}} \leq d(v,X_j) = h(v,X_j)\) by property~\((p2).\) Thus the weight of~\((X_j,v_{min})\) is also less than that of~\((v,X_j)\) and so~(\ref{mst_prop}) is true for~\(e=(v,X_j).\)~\(\qed\)

\emph{Proof of~(\ref{add_one_est1}) in Lemma~\ref{var_low_lem1}}: From~(\ref{add_one}), we get that if the event mentioned in the statement of the Lemma occurs, then the weight~\[W({\cal T}_{n+1}) = MST_{n+1}  = MST_{n}(j) + h^{\alpha}(X_j,v_{min}) = MST_{n}(j) + d^{\alpha}(X_j,v_{min}) \] since~\(X_j \in S_0.\) The vertex~\(v_{min}\) is necessarily one of the eight vertices in\\\(\{v_i\}_{1 \leq i \leq 12} \setminus \{v_1,v_4,v_7,v_{10}\}.\) By construction, the distance between any node within the square~\(R_i\) (labelled~\(E\) in Figure~\ref{good_square}) and any node within the square labelled~\(3,\) for example, is at least~\((3g-1)\frac{A}{\sqrt{n}}\) and at most~\((5g-1)\frac{A}{\sqrt{n}}.\)~\(\qed\)

%NEED TO IMPOSE CONDITIONS ON H SO THAT THE WEIGHTS ARE DIFF?? AND TO SEE CRFLLY+ ETC..

\subsection*{\em Proof of Lemma~\ref{var_egj_lem}}
For a constant~\(\theta > 0\) let~\(E_j(g,\theta)\) be the event that the number of~\((g,j)-\)good squares in~\(\{R_k\}\) is at least~\(\theta \cdot n.\) Recalling the node distribution parameters~\(\epsilon_1\) and~\(\epsilon_2\) from~(\ref{f_eq}), we first show that there are positive constants~\(\theta_1 = \theta_1(g,s_0,\epsilon_1,\epsilon_2)\) and~\(\theta_2 = \theta_2(g,s_0,\epsilon_1,\epsilon_2)\) such that for any~\(1 \leq j \leq n+1,\)
\begin{equation}\label{egj_est}
\mathbb{P}\left(E_j(g,\theta_1)\right) \geq 1-e^{-\theta_2 n}.
\end{equation}
We then let
\begin{equation}\label{thet_not_def}
\theta_0 := \min(\theta_1(2g,s_0,\epsilon_1,\epsilon_2),\theta_1(g,s_0,\epsilon_1,\epsilon_2))
\end{equation}
and define \[E_{tot}(j) := E_{j}(2g,\theta_0) \cap F_j(2g) \cap E_j(g,\theta_0) \cap F'_j(g),\] where we recall that~\(F_j(2g)\) denotes the event that the node~\(X_j\) belongs to a~\((2g,j)-\)good square in~\(\{R_i\}\) and an analogous definition holds for~\(F'_j(g)\) with  the node~\(X'_j.\) We show that there exists a constant~\(D = D(g,s_0,\epsilon_1,\epsilon_2) > 0\) such that for any~\(1 \leq j \leq n+1\)
\begin{equation}\label{ej_fj_est}
\mathbb{P}(E_{tot}(j)) \geq D,
\end{equation}
completing the proof of Lemma~\ref{var_egj_lem}.

We prove~(\ref{egj_est}) and~(\ref{ej_fj_est}) in that order below.\\
\emph{Proof of~(\ref{egj_est})}: We use Poissonization and let~\({\cal P}\) be a Poisson process of intensity~\(nf(.)\) in the unit square, with corresponding probability measure~\(\mathbb{P}_0.\) Analogous to the definition of~\((g,j)-\)good square, we define~\(R_i\) to be~\(g-\)good square if the twelve~\(\frac{A}{\sqrt{n}} \times \frac{A}{\sqrt{n}}\) squares as in Figure~\ref{good_square}\((a)\) each contain exactly one node of~\({\cal P}\) and the rest of big square of size~\(\frac{(30g+1)A}{\sqrt{n}} \times \frac{(30g+1)A}{\sqrt{n}}\) is empty.
The number of nodes~\(N(R_i)\) of~\({\cal P}\) in the square~\(R_i\) is Poisson distributed with mean~\(n\int_{R_i}f(x) dx \in [\epsilon_1 A^2, \epsilon_2 A^2]\) (see the bounds for~\(f(.)\) in~(\ref{f_eq})) and so~\(\mathbb{P}_0(R_i \text{ is }g-\text{good}) \geq p_0\) for some constant~\(p_0  = p_0(\epsilon_1,\epsilon_2,g) > 0\) not depending on~\(i.\)

If~\(R_i\) and~\(R_k\) are two squares such that the corresponding neighbourhoods~\({\cal N}_{15g}(R_i) \cap {\cal N}_{15g}(R_k) = \emptyset,\) then the events~\(\{R_i \text{ is } g-\text{good}\}\) and~\(\{R_k \text{ is } g-\text{good}\}\) are independent since the Poisson process is independent on disjoint sets. We therefore pick a maximal set of squares~\(\{R_l\}_{l \in {\cal Q}}\) whose\\\(15g-\)neighbourhoods are empty. Since the area of the square~\(S_0\) is~\(s_0^2,\) a constant and~\({\cal N}_{15}(R_i)\) contains~\((30g+1)^2\) squares in~\(\{R_l\},\) each of area~\(\frac{A^2}{n},\) we have that~\(\#{\cal Q} \geq \frac{1}{2}\cdot \frac{s_0^2n}{A^2(30g+1)^2} \geq \frac{s_0^2n}{8(30g+1)^2}\) using~\(A \leq 1+\frac{1}{\log{n}} \leq 2.\)

If~\(N_{good} := \sum_{i \in {\cal Q}} \ind(R_i \text{ is } g-\text{good})\) denotes the number of~\(g-\)good squares in the collection~\({\cal Q},\)
then~\(N_{good}\) is a sum of independent Bernoulli random variables, each with mean at least~\(p_0.\) Therefore~\(\mathbb{E}_0(N_{good}) \geq \frac{s_0^2p_0n}{8(30g+1)^2} =: 2\theta_1 n\) and using the standard deviation estimate~(\ref{std_dev_down}) (see Appendix) with~\(m = \#{\cal Q}, \mu_1 = p_0\) and~\(\epsilon = \frac{1}{2},\) we get \[\mathbb{P}_0\left(N_{good} \geq \theta_1 n\right) \geq 1-e^{-2\theta_2 n}\] for some positive constant~\(\theta_2 = \theta_2(g,s_0,\epsilon_1,\epsilon_2).\) Using the dePoissonization formula~(\ref{de_poiss_ax}), there exists a constant~\(D > 0\) not depending on~\(j\) such that~\[\mathbb{P}(E_j(g,\theta_1)) \geq 1-D\sqrt{n} \cdot e^{-2\theta_2 n} \geq 1-e^{-\theta_2n}\] for all~\(n\) large.~\(\qed\)

\emph{Proof of~(\ref{ej_fj_est})}: We let~\({\cal S}_{good}(2g)\) and~\({\cal S}_{good}(g)\) respectively denote the random set of all~\((2g,j)-\)good squares
and the random set of all~\((g,j)-\)good squares in~\(\{R_i\}.\) If~\(E_j(2g,\theta_0) \cap E_j(g,\theta_0)\) occurs where~\(\theta_0\) is as defined prior to~(\ref{ej_fj_est}), then~\(\#{\cal S}_{good}(2g) \geq \theta_0\cdot n\)
and~\(\#{\cal S}_{good}(g) \geq \theta_0\cdot n.\) The event~\(E_j(2g,\theta_0) \cap E_j(g,\theta_0)\) is independent of~\(X_j\) and~\(X'_j\)
and so
\begin{eqnarray}
\mathbb{P}(E_{tot}(j)) &=& \sum_{{\cal S}_1 : \#{\cal S}_1 \geq \theta_0 n} \sum_{{\cal S}_2 : \#{\cal S}_2 \geq \theta_0n} \mathbb{P}\left(\{X_j \in {\cal S}_1\} \cap \{X'_j \in {\cal S}_2\}  \right. \nonumber\\
&&\;\;\;\;\;\left. \cap \{{\cal S}_{good}(2g) = {\cal S}_1\} \cap \{{\cal S}_{good}(g) = {\cal S}_2\}\right) \nonumber\\
&=& \sum_{{\cal S}_1 : \#{\cal S}_1 \geq \theta_0 n} \sum_{{\cal S}_2 : \#{\cal S}_2 \geq \theta_0n} \mathbb{P}(X_j \in {\cal S}_1) \mathbb{P}(X'_j \in {\cal S}_2) \nonumber\\
&&\;\;\;\;\; \mathbb{P}\left(\{{\cal S}_{good}(2g) = {\cal S}_1\} \cap \{{\cal S}_{good}(g) = {\cal S}_2\}\right). \label{e_tot_temp}
\end{eqnarray}

For any collection~\({\cal S}_1\) containing~\({\cal S}_1 \geq \theta_0n\) squares from~\(\{R_i\},\) we have that~\[\mathbb{P}(X_j \in {\cal S}_1) \geq \int_{{\cal S}_1} f(x)dx \geq \epsilon_1 \cdot\frac{A^2}{n} \cdot \theta_0n \geq \epsilon_1 \theta_0,\] since~\(A \geq 1.\) Similarly~\(\mathbb{P}(X'_j \in {\cal S}_2) \geq \epsilon_1 \theta_0\) and so from~(\ref{e_tot_temp}) we get that
\begin{eqnarray}
\mathbb{P}(E_{tot}(j)) &\geq& (\epsilon_1 \theta_0)^2 \sum_{{\cal S}_1 : \#{\cal S}_1 \geq \theta_0 n} \sum_{{\cal S}_2 : \#{\cal S}_2 \geq \theta_0n} \mathbb{P}\left(\{{\cal S}_{good}(2g) = {\cal S}_1\} \cap \{{\cal S}_{good}(g) = {\cal S}_2\}\right) \nonumber\\
&=& (\epsilon_1 \theta_0)^2 \mathbb{P}\left(E_j(2g,\theta_0) \cap E_j(g,\theta_0)\right) \nonumber\\
&\geq& (\epsilon_1 \theta_0)^2 (1-e^{-n\cdot \theta_2(2g,s_0,\epsilon_1,\epsilon_2)}-e^{-n\cdot \theta_2(g,s_0,\epsilon_1,\epsilon_2)}) \nonumber
\end{eqnarray}
using~(\ref{egj_est}) and the definition of~\(\theta_0\) in~(\ref{thet_not_def}).~\(\qed\)
%and some constant~\(D > 0\) not depending on~\(j.\)

%TO SEE CRFLY +eTC...

\setcounter{equation}{0}
\renewcommand\theequation{\thesection.\arabic{equation}}
\section{Convergence properties for MST}\label{mst_conv}
In this section, we study convergence properties for the MST weight~\(MST_n\) as defined in~(\ref{min_weight_tree}), appropriately scaled and centred. The following is the main result of this section.
\begin{Theorem}\label{thm_mst_conv} For every~\(\alpha > 0\) we have that~\[\frac{1}{n^{1-\frac{\alpha}{2}}}(MST_n-\mathbb{E}MST_n) \longrightarrow 0 \text{ a.s.\ }\] as~\(n \rightarrow \infty.\)
\end{Theorem}
The proof is standard and follows from subsequence arguments (Steele (1988)). For completeness we provide a
proof below.

\emph{Proof of Theorem~\ref{thm_mst_conv}}: We prove the almost sure convergence via a subsequence argument
using the variance upper bound for~\(MST_n\) obtained in Theorem~\ref{thm_mst_var_up}. Indeed, from Theorem~\ref{thm_mst_var_up} we have that~\[var\left(\frac{MST_n}{n^{1-\frac{\alpha}{2}}}\right) \leq \frac{C}{n}\]
and so an application of the Borel-Cantelli Lemma gives us that~\[\frac{1}{n^{2-\alpha}}(MST_{n^2} - \mathbb{E}MST_{n^2}) \longrightarrow 0 \text{ a.s.\ } \] as~\(n \rightarrow \infty.\) To prove convergence along the subsequence~\(a_n = n,\) we let
\begin{equation}\label{d_def_mst}
D_n := \max_{n^2 \leq k < (n+1)^2} \left|MST_k - MST_{n^2}\right|
\end{equation}
and show that~\(\mathbb{E} D^2_n \leq C n^{2-2\alpha}\) for some constant~\(C > 0.\) This would then imply that~\(\left(\frac{\mathbb{E}D_n}{n^{2-\alpha}}\right)^2 \leq \frac{\mathbb{E}D^2_n}{n^{4-2\alpha}} \leq \frac{C}{n^{2}} \longrightarrow 0\) as~\(n \rightarrow \infty\)
and we also get from Borel-Cantelli Lemma that~\(\frac{D_n}{n^{2-\alpha}} \longrightarrow 0\) a.s.\
as~\(n \rightarrow \infty.\) For~\(n^2 \leq k < (n+1)^2\) we then write
\begin{eqnarray}
\frac{|MST_k - \mathbb{E}MST_{k}|}{k^{1-\frac{\alpha}{2}}} &\leq& \frac{|MST_k - MST_{n^2}|}{k^{1-\frac{\alpha}{2}}} + \frac{\mathbb{E}|MST_{k} - MST_{n^2}|}{k^{1-\frac{\alpha}{2}}} \nonumber\\
&\leq& \frac{D_n}{k^{1-\frac{\alpha}{2}}} + \frac{\mathbb{E}D_n}{k^{1-\frac{\alpha}{2}}} \nonumber\\
&\leq& \frac{D_n}{n^{2-\alpha}} + \frac{\mathbb{E}D_n}{n^{2-\alpha}} \nonumber
\end{eqnarray}
and get that~\(\frac{MST_{k} - \mathbb{E}MST_k}{k^{1-\frac{\alpha}{2}}} \longrightarrow 0\) a.s.\ as~\(n \rightarrow \infty.\)

To estimate~\(D_n\) we use the one node difference estimate~(\ref{cruc_est}) to get that~\(|MST_{k+1} - MST_k| \leq f_{1,k} + f_{2,k}\)
where~\(f_{1,k}\) and~\(f_{2,k}\) are such that
\begin{equation}\label{ef_est}
\left(\mathbb{E}f_{1,k}\right)^2 \leq \mathbb{E}f^2_{1,k} \leq \frac{C}{k^{\alpha}} \text{ and } \left(\mathbb{E}f_{2,k}\right)^2 \leq \mathbb{E}f^2_{2,k} \leq \frac{C}{k^{\alpha}}
\end{equation}
for some constant~\(C > 0\) (see~(\ref{f1_prop1}) and~(\ref{f2_prop1})). Thus telescoping we get
\begin{equation} \nonumber
|MST_{k}-MST_{n^2}| \leq \sum_{l=n^2}^{k} f_{1,l} + f_{2,l}
\end{equation}
and so \[D_n = \max_{n^2 \leq k <(n+1)^2} |MST_k - MST_{n^2}| \leq \sum_{l=n^2}^{(n+1)^2} f_{1,l} + f_{2,l}.\]

Using~\((\sum_{i=1}^{t} a_i)^2 \leq t \sum_{i=1}^{t} a_i^2\) we get that~\(\mathbb{E}D_n^2 \) is bounded above by
\[((n+1)^2-n^2)\sum_{l=n^2}^{(n+1)^2}\mathbb{E}(f_{1,l} + f_{2,l})^2 \leq  2((n+1)^2-n^2)\sum_{l=n^2}^{(n+1)^2}(\mathbb{E}f^2_{1,l} + \mathbb{E}f^2_{2,l})\] and plugging the estimates from~(\ref{ef_est}) we finally get
\begin{equation}\label{fn_est1}
\mathbb{E} D_n^2 \leq 2((n+1)^2-n^2)\sum_{l=n^2}^{(n+1)^2} \frac{2C}{l^{\alpha}} \leq 2((n+1)^2-n^2)^2\frac{2C}{n^{2\alpha}} \leq 8C n^{2-2\alpha},
\end{equation}
proving the desired estimate for~\(\mathbb{E}D_n^2.\)~\(\qed\)

\setcounter{equation}{0}
\renewcommand\theequation{\thesection.\arabic{equation}}
\section{Uniform MSTs}\label{mst_uni}
In this Section, we assume that the nodes~\(\{X_i\}_{1 \leq i \leq n}\) are uniformly distributed in the unit square and obtain bounds on the asymptotic values of the expected weight, appropriately scaled and centred. We assume that the positive edge weight function~\(h : \mathbb{R}^2 \times \mathbb{R}^2 \rightarrow \mathbb{R}\) satisfies~(\ref{eq_met1}) along with the following two properties:\\
\((b1)\) For every~\(a > 0\) we have
\[h(au,av) = a\cdot h(u,v) \text{ for all } u,v \in \mathbb{R}^2\]
\((b2)\) There exists a constant~\(h_0 > 0\) such that for all~\(b \in \mathbb{R}^2\) we have
\begin{equation}\label{h0_def}
h(b+u,b+v) \leq h_0 \cdot h(u,v).
\end{equation}
For example, recalling that~\(d(u,v)\) denotes the Euclidean distance between~\(u\) and~\(v,\) we have that the function
\begin{equation}\label{h_ex}
h(u,v) = d(u,v) + \frac{1}{2} |d(u,0) - d(v,0)|
\end{equation}
is a metric since
\[d(u,0) \leq d(0,v) + d(v,u) \text{ and } d(v,0) \leq d(0,u) + d(u,v)\] by triangle inequality and so
\[d(u,v) \leq h(u,v) \leq d(u,v) + \frac{1}{2} d(u,v).\] This implies that~\(h(u,v)\) satisfies~(\ref{eq_met1}) and by definition~\(h\) also
satisfies~\((b1).\) Moreover using the triangle inequality we have
\begin{eqnarray}
h(b+u,b+v) &=& d(b+u,b+v) + \frac{1}{2}|d(u+b,0)-d(v+b,0)|  \nonumber\\
&=& d(u,v) + \frac{1}{2}|d(u+b,0)-d(v+b,0)| \nonumber\\
&\leq& d(u,v) + \frac{1}{2} d(u+b,v+b) \nonumber\\
&=&d(u,v) + \frac{1}{2} d(u,v) \nonumber\\
&\leq& \frac{3}{2} h(u,v) \nonumber
\end{eqnarray}
by definition of~\(h\) in~(\ref{h_ex}). Thus~\((b2)\) is also satisfied with~\(h_0 = \frac{3}{2}.\)

We have the following result.
\begin{Theorem}\label{mst_unif}
Suppose the distribution function~\(f(.)\) is uniform, i.e.,~\(\epsilon_1 = \epsilon_2 = 1\) in~(\ref{f_eq}) and the edge weight function~\(h(u,v)\) satisfies~(\ref{eq_met1}) and properties~\((b1)-(b2)\) above.  We then have
\begin{equation}\label{lim_sup_est}
0 < \liminf_{n} \frac{\mathbb{E} MST_{n}}{n^{1-\frac{\alpha}{2}}} \leq \limsup_{n} \frac{\mathbb{E} MST_{n}}{n^{1-\frac{\alpha}{2}}} \leq h_0^{\alpha} \cdot \liminf_{n} \frac{\mathbb{E} MST_{n}}{n^{1-\frac{\alpha}{2}}} < \infty.
\end{equation}
\end{Theorem}
Thus the scaled weight of the minimum spanning tree remains bounded within a factor of~\(h_0^{\alpha}.\)

We begin with the following Lemma.
\begin{Lemma}\label{mst_gen_lem} There is a constant~\(D > 0\) such that for any positive integers~\(n_1,n_2 \geq 1\) we have
\begin{eqnarray}
\mathbb{E} MST_{n_1+n_2} \leq \mathbb{E}MST_{n_1} + n_2\left(\frac{D}{n_1}\right)^{\frac{\alpha}{2}}. \label{mst_ab2}
\end{eqnarray}
Moreover, for any fixed integer~\(m \geq 1\) we have
\begin{equation}\label{km_su}
\limsup_n \frac{\mathbb{E} MST_{n}}{n^{1-\frac{\alpha}{2}}} \leq \limsup_k \frac{\mathbb{E}MST_{km}}{(km)^{1-\frac{\alpha}{2}}}.
\end{equation}
\end{Lemma}

\emph{Proof of~(\ref{mst_ab2}) in Lemma~\ref{mst_gen_lem}}: Let~\({\cal T}_1\) be the MST formed by the~\(n_1\) nodes~\(\{X_i\}_{1 \leq i \leq n_1}.\)
We join each~\(X_i, n_1+1 \leq i \leq n_1+n_2\) to the node in~\(\{X_j\}_{1 \leq j \leq n_1}\) closest to~\(X_i\) in Euclidean distance
by an edge~\(e_i,\) whose length is~\(d^{\alpha}(X_i,\{X_j\}_{1 \leq j \leq n_1}),\)
where we recall that~\(d(X_i,\{X_j\}_{1 \leq j \leq n_1})\) is the minimum distance
between~\(X_i\) and the nodes in~\(\{X_j\}_{1 \leq j \leq n_1}.\) From~(\ref{eq_met1}),
we have that the weight of~\(e_i\) is at most~\(c_2^{\alpha}d^{\alpha}(X_i,\{X_j\}_{1 \leq j \leq n_1}).\)

The tree~\({\cal T}_1 \cup \{e_i\}_{n_1+1 \leq i \leq n_1+n_2}\) contains all the~\(n_1+n_2\) nodes and so
we have
\begin{eqnarray}
\mathbb{E}MST_{n_1+n_2} &\leq& \mathbb{E}MST_{n_1} + c_2^{\alpha}\sum_{i=n_1+1}^{n_1+n_2} \mathbb{E}d^{\alpha}(X_i,\{X_j\}_{1 \leq j \leq n_1}) \nonumber\\
&=& \mathbb{E}MST_{n_1} + n_2 \cdot c_2^{\alpha} \cdot \mathbb{E}d^{\alpha}(X_{n_1+1},\{X_j\}_{1 \leq j \leq n_1}). \nonumber
\end{eqnarray}
From the estimate~(\ref{f1_prop1}), we have that~\(\mathbb{E}d^{\alpha}(X_{n_1+1},\{X_j\}_{1 \leq j \leq n_1}) \leq \left(\frac{D_1}{n_1}\right)^{\frac{\alpha}{2}}\) for some constant~\(D_1 > 0\) not depending on~\(n_1\) or~\(n_2\) and this proves~(\ref{mst_ab2}).~\(\qed\)

\emph{Proof of~(\ref{km_su}) in Lemma~\ref{mst_gen_lem}}: Fix an integer~\(m \geq 1\) and write~\(n = qm + s\)
where~\(q = q(n) \geq 1\) and~\(0 \leq s = s(n) \leq m-1\) are integers.
As~\(n \rightarrow \infty,\)
\begin{equation}\label{kn_def}
q(n) \longrightarrow \infty \text{ and }\frac{n}{q(n)} \longrightarrow m.
\end{equation}

Using property~\((t1)\) with~\(n_1 = qm\) and~\(n_2 = s\) to get that
\[\mathbb{E}MST_{n} = \mathbb{E}MST_{qm+s} \leq \mathbb{E}MST_{qm} + s\left(\frac{D}{qm}\right)^{\frac{\alpha}{2}} \leq \mathbb{E}MST_{qm} + m\left(\frac{D}{qm}\right)^{\frac{\alpha}{2}},\]
since~\(s < m\) and so
\begin{equation}
\limsup_n \frac{\mathbb{E}MST_{n}}{n^{1-\frac{\alpha}{2}}} \leq \limsup_n \left(\frac{qm}{n}\right)^{1-\frac{\alpha}{2}} \frac{\mathbb{E} MST_{qm}}{(qm)^{1-\frac{\alpha}{2}}}+\limsup_n \frac{m}{n^{1-\frac{\alpha}{2}}} \left(\frac{D}{qm}\right)^{\frac{\alpha}{2}}.
\label{mst_temp_km1}
\end{equation}
Since~\(\frac{n}{q(n)} \longrightarrow m, q(n) \longrightarrow \infty\) as~\(n \rightarrow \infty,\) (see~(\ref{kn_def})) and~\(m\) is fixed, we have that~\[\frac{m}{n^{1-\frac{\alpha}{2}}} \left(\frac{D}{qm}\right)^{\frac{\alpha}{2}} = \frac{m^{1-\frac{\alpha}{2}}D^{\frac{\alpha}{2}}}{n}\cdot \left(\frac{n}{q}\right)^{\frac{\alpha}{2}} \longrightarrow 0\] as~\(n \rightarrow \infty.\)
Thus the second term in the right side of~(\ref{mst_temp_km1}) is zero
and the first term in the right side of~(\ref{mst_temp_km1}) equals~\(\limsup_n \frac{\mathbb{E} MST_{q m}}{(qm)^{1-\frac{\alpha}{2}}}.\)

But since~\(q(n) \geq \frac{n-m}{m} \geq \frac{l-m}{m}\) for~\(n \geq l,\) we have that
\begin{equation}\label{gtu}
\sup_{n \geq l} \frac{\mathbb{E} MST_{qm}}{(qm)^{1-\frac{\alpha}{2}}} = \sup_{n \geq l} \frac{\mathbb{E} MST_{q(n)m}}{(q(n)m)^{1-\frac{\alpha}{2}}} \leq \sup_{k \geq \frac{l-m}{m}} \frac{\mathbb{E} MST_{km}}{(km)^{1-\frac{\alpha}{2}}}
\end{equation}
and as~\(l \uparrow \infty,\) the final term in~(\ref{gtu}) converges to the second term in~(\ref{km_su}).~\(\qed\)

%In particular if~\(l \geq 1\) is an integer then \[\sup_{ n \geq lm + m} \frac{\mathbb{E} MST_{q^2 m^2}}{qm} \geq L_1.\]

If~\(\lambda := \liminf_n \frac{\mathbb{E}MST_{n}}{n^{1-\frac{\alpha}{2}}}\) then from~(\ref{exp_mst_bound}), we have that~\(\lambda > 0.\) We show below that for every~\(\epsilon > 0,\) there exists~\(m\) sufficiently large such that
\begin{equation}\label{second_step_est}
\limsup_{k} \frac{\mathbb{E} MST_{km}}{(km)^{1-\frac{\alpha}{2}}} \leq h_0^{\alpha} \cdot \lambda + \epsilon,
\end{equation}
where~\(h_0\) is as in~(\ref{h0_def}). Since~\(\epsilon > 0\) is arbitrary, Theorem~\ref{mst_unif} then follows from~(\ref{second_step_est}) and~(\ref{km_su}). To prove~(\ref{second_step_est}), let~\(k\) and~\(m\) be positive integers and distribute~\(km\) nodes~\(\{X_i\}_{1 \leq i \leq km}\) independently
and uniformly in the unit square~\(S.\) Also divide~\(S\) into~\(k\) disjoint
squares~\(\{W_j\}_{1 \leq j \leq \frac{k}{A_k^2}}\) each of size~\(\frac{A_k}{\sqrt{k}} \times \frac{A_k}{\sqrt{k}}\)
where~\(A_k \in \left[1, 1 + \frac{1}{\log{k}}\right]\) is such that~\(\frac{\sqrt{k}}{A_k}\)
is an integer for all large~\(k \geq K_0,\) not depending on~\(m.\) This is possible by an analogous argument as in~(\ref{poss}).
Further, we label the squares~\(\{W_j\}\) as in Figure~\ref{fig_squares}
so that the top left most square is labelled~\(W_1,\) the square below~\(W_1\) is~\(W_2\) and so on. %until
%we reach the square~\(W_{\frac{\sqrt{k}}{A_k}}\) intersecting the bottom edge of the unit square~\(S.\) The square to the right of~\(W_k\)
%is then labelled~\(W_{k+1}\) and the square above~\(W_{k+1}\) is~\(W_{k+2}\) and so on.

For~\(1 \leq j \leq \frac{k}{A_k^2},\) let~\(N(j)\) be the number of nodes of~\(\{X_i\}_{1 \leq i \leq km}\)
present in the square~\(W_j\) and let~\({\cal T}(N(j))\) be the MST containing all the~\(N(j)\) nodes
with corresponding length~\(MST_j(N(j)).\) Also, for~\(1 \leq j \leq \frac{k}{A_k^2}-1,\)
we let~\(j~+~T^{next}_j := \min\{i \geq j+1 : W_i \text{ is not empty}\}\) be the next nonempty square after~\(W_j\)
and let~\(e_j\) be the edge with one endvertex in~\(W_j\) and the other endvertex in~\(W_{j+T^{next}_j},\)
having the smallest Euclidean length. We denote the Euclidean length of~\(e_j\) to be~\(d(e_j)\)
and set~\(MST_j(N(j)) = d(e_j) = 0\) if~\(W_j\) is empty.  The union~\(\cup_{1 \leq j \leq \frac{k}{A_k^2}} {\cal T}(j) \cup \cup_{1 \leq j \leq \frac{k}{A_k^2}-1} \{e_j\}\) is a spanning tree containing all the~\(km\) nodes and so
\begin{eqnarray}\label{mst_km}
MST_{km} \leq \sum_{j=1}^{\frac{k}{A_k^2}} MST_j(N(j))  + c_2^{\alpha} \sum_{j=1}^{\frac{k}{A_k^2}-1} d^{\alpha}(e_j),
\end{eqnarray}
using~(\ref{eq_met1}). In Appendix, we use the translation property~\((b2)\) to get that
\[MST_j(N(j)) \leq h_0^{\alpha} \cdot MST(N(j)),\]
where~\(MST(N(j))\) is the MST weight of the configuration of nodes present in~\(W_j\) with the centre of the square~\(W_j\)
shifted to the origin: i.e., if~\(u_1,\ldots,u_w\) are the nodes of~\(\{X_i\}\) present in the square~\(W_j\) whose centre is~\(s_j,\)
then~\(MST(N(j))\) is the MST weight of complete graph formed by the nodes~\(u_1-s_j,\ldots,u_w-s_j.\)

Taking expectations we get
\begin{eqnarray}\label{mst_km2}
\mathbb{E}MST_{km} &\leq& h_0^{\alpha}\sum_{j=1}^{\frac{k}{A_k^2}} \mathbb{E}MST(N(j)) + c_2^{\alpha} \cdot \frac{k}{A_k^2} \max_{1 \leq j \leq \frac{k}{A_k^2}-1} \mathbb{E}d^{\alpha}(e_j) \nonumber\\
&=& h_0^{\alpha} \cdot \frac{k}{A_k^2} \mathbb{E}MST(N(1)) + c_2^{\alpha} \cdot \frac{k}{A_k^2} \max_{1 \leq j \leq \frac{k}{A_k^2}-1} \mathbb{E}d^{\alpha}(e_j)
\end{eqnarray}
and for convenience we write~\[\mathbb{E}MST(N(1)) = I_1 + I_2,\] where~\[I_1 =\mathbb{E} MST(N(1))\ind(F_1),I_2 = \mathbb{E}MST(N(1))\ind(F_1^c)\] and~\[F_1 := \{mA_k^2 - \sqrt{m}\log{m} \leq N_1 \leq mA_k^2 + \sqrt{m}\log{m}\}\] to get
\begin{eqnarray}\label{mst_km3}
\mathbb{E}MST_{km} \leq h_0^{\alpha} \cdot \frac{k}{A_k^2} (I_1+I_2)  + c_2^{\alpha} \cdot \frac{k}{A_k^2} \max_{1 \leq j \leq \frac{k}{A_k^2}-1} \mathbb{E}d^{\alpha}(e_j)
\end{eqnarray}

The following Lemma estimates each sum in~(\ref{mst_km3}) starting with the second term.
\begin{Lemma}\label{i1_i2_lem}
There are positive constants~\(D_1,D_2\) not depending on~\(k\) or~\(m\) such that
\begin{equation}\label{s_eval}
\max_{1 \leq j \leq \frac{k}{A_k^2}-1} \mathbb{E}d^{\alpha}(e_j) \leq D_1\left(\frac{A_k}{\sqrt{k}}\right)^{\alpha}\left(\left(\frac{\log{m}}{\sqrt{m}}\right)^{\alpha} + e^{-2\sqrt{m}\log{m}}\right),
\end{equation}
\begin{equation}
I_1 \leq \left(\frac{A_k}{\sqrt{k}}\right)^{\alpha} \left(\mathbb{E} MST_{m} + \frac{D_1 \cdot m^{1-\frac{\alpha}{2}}}{\log{k}} + D_1\cdot  m^{\frac{1-\alpha}{2}} (\log{m})\right) \label{i1_est}
\end{equation}
and
\begin{equation}\label{i2_est}
I_2 \leq D_1 \left(\frac{A_k}{\sqrt{k}}\right)^{\alpha}\left(e^{-m \cdot D_2} + m^{1-\frac{\alpha}{2}} \cdot \frac{1}{(\log{m})^2}\right).
\end{equation}
\end{Lemma}

Substituting the estimates~(\ref{s_eval}),~(\ref{i1_est}) and~(\ref{i2_est}) of Lemma~\ref{i1_i2_lem} into~(\ref{mst_km3}), we get that~\[\mathbb{E}MST_{km} \leq k^{1-\frac{\alpha}{2}} A_k^{\alpha-2}\left( h_0^{\alpha} \cdot \mathbb{E} MST_{m} + D_1 \cdot R_m(k)\right)\]  where~\(R_m(k)\) equals \[\left(\frac{m^{1-\frac{\alpha}{2}}}{\log{k}} + m^{\frac{1-\alpha}{2}} (\log{m})\right) + \left(\frac{m^{1-\frac{\alpha}{2}} }{(\log{m})^2} + e^{-D_2 \cdot m}\right) +\left(\left(\frac{\log{m}}{\sqrt{m}}\right)^{\alpha} + e^{-\sqrt{m} \log{m}}\right).\]
Using~\(1 \leq A_k \leq 1+\frac{1}{\log{k}} \longrightarrow 1\) as~\(k \rightarrow \infty\)
and absorbing~\(e^{-D_2 \cdot m}\) into~\(e^{-\sqrt{m} \log{m}},\) we then get
that~
\begin{equation}\label{eq_no}
\limsup_k \frac{\mathbb{E} MST_{km}}{(km)^{1-\frac{\alpha}{2}}} \leq h_0^{\alpha} \cdot \frac{\mathbb{E} MST_{m}}{m^{1-\frac{\alpha}{2}}} + 2D_1\cdot R'_m
\end{equation}
where
\begin{equation}\nonumber
R'_m := \frac{\log{m}}{\sqrt{m}} + \frac{1}{(\log{m})^2}  + \frac{(\log{m})^{\alpha}}{m} + m^{\frac{\alpha}{2}-1}e^{-\sqrt{m}\log{m}} \leq \epsilon
\end{equation}
for any~\(\epsilon > 0\) and all~\(m\) large.  Letting~\(\{m_j\}\) be any sequence such that\\\(\frac{\mathbb{E} MST_{m_j}}{m_j^{1-\frac{\alpha}{2}}} \longrightarrow \lambda = \liminf_n \frac{\mathbb{E}MST_n}{n^{1-\frac{\alpha}{2}}}\) and allowing~\(m \rightarrow \infty\) through the sequence~\(\{m_j\}\) we get from~(\ref{km_su}), ~(\ref{eq_no}) and the above discussion that~\[\limsup_n \frac{\mathbb{E}MST_n}{n^{1-\frac{\alpha}{2}}} \leq h_0^{\alpha} \cdot \lambda + 2D_1\epsilon.\] Since~\(\epsilon  >0\) is arbitrary, we get~(\ref{second_step_est}).\(\qed\)

\begin{figure}[tbp]
\centering
%\fbox{
\includegraphics[width=3in, trim= 60 300 80 250, clip=true]{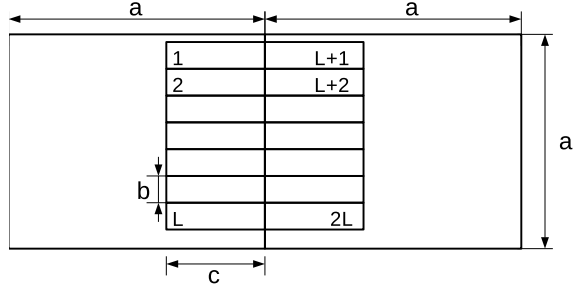}
%}
\caption{Determining the length of the edge~\(e_j\) with one endvertex in the left~\(a\times a\) square~\(W_j\) and other endvertex in right~\(a\times a\) square~\(W_{j+1}.\) Here~\(a = \frac{A_k}{\sqrt{k}}, b = \frac{A_k}{\sqrt{km}}\) and~\(c= \frac{4A_k\log{m}}{\sqrt{km}}.\)}
\label{fig_w_rect}
\end{figure}

\emph{Proof of~(\ref{s_eval}) in Lemma~\ref{i1_i2_lem}}: We let~\(1 \leq j \leq \frac{k}{A_k^2}-1\) and suppose that the right edge of~\(W_j\) is the left edge of~\(W_{j+1}.\) Let~\(L\) be the largest integer less than or equal to~\(\sqrt{m}\) and for~\(1 \leq i \leq L\) let~\(W^{small}_j(i)\) and~\(W^{small}_{j+1}(i)\) be disjoint~\(\frac{4A_k\log{m}}{\sqrt{km}} \times \frac{A_k}{\sqrt{km}}\) rectangles contained in~\(W_j\) and~\(W_{j+1}\) respectively,
and sharing a common edge as shown in Figure~\ref{fig_w_rect}. Here the two~\(a \times a\) squares represent~\(W_j\) and~\(W_{j+1}\) with~\(a = \frac{A_k}{\sqrt{k}}, b = \frac{A_k}{\sqrt{km}}\) and~\(c= \frac{4A_k\log{m}}{\sqrt{km}}.\) The~\(c \times b\) rectangle labelled~\(i\) represents~\(W^{small}_j(i)\) for~\(1 \leq  i \leq L\) and represents~\(W^{small}_{j+1}(i-L)\) for~\(L+1 \leq i \leq 2L.\)

Let~\(E_{nice}(j)\) be the event that there exists a pair of rectangles~\(W^{small}_j(i)\) and~\(W^{small}_{j+1}(i)\) both of which contain
at least one of the nodes of~\(\{X_u\}_{1 \leq u \leq km}.\) If~\(E_{nice}(j)\) occurs, then the minimum length of an edge
with one endvertex in~\(W_j\) and other endvertex in~\(W_{j+1}\) is at most~\(2c+2b \leq 4c = \frac{16A_k\log{m}}{\sqrt{km}}.\)
If~\(E_{nice}(j)\) does not occur, we use the upper bound~\(d(e_j) \leq \frac{2A_kT^{next}_j}{\sqrt{k}}\) and so
\[d(e_j) \leq \frac{16A_k\log{m}}{\sqrt{km}}\ind(E_{nice}(j)) + \frac{2A_kT^{next}_j}{\sqrt{k}} \ind(E^c_{nice}(j)).\] Using~\((a+b)^{\alpha} \leq 2^{\alpha}(a^{\alpha}+ b^{\alpha})\) for all~\(a,b,\alpha > 0\) we therefore get that
\begin{equation}
\mathbb{E}d^{\alpha}(e_j) \leq \left(\frac{4A_k}{\sqrt{k}}\right)^{\alpha}\left(\left(\frac{8\log{m}}{\sqrt{m}}\right)^{\alpha}\mathbb{P}(E_{nice}(j)) + \mathbb{E}\left(T^{next}_j\right)^{\alpha}\ind(E^c_{nice}(j))\right).\label{ed_alpha}
\end{equation}

To estimate the second term within the brackets in~(\ref{ed_alpha}), we first use Cauchy-Schwarz inequality and write~
\begin{equation}\label{ed_alp2}
\mathbb{E}\left(T^{next}_j\right)^{\alpha}\ind(E^c_{nice}(j)) \leq \left(\mathbb{E}\left(T^{next}_j\right)^{2\alpha}\right)^{\frac{1}{2}} \left(\mathbb{P}(E_{nice}^c(j))\right)^{\frac{1}{2}}.
\end{equation}
For any~\(l \geq 1\) we have~\(T^{next}_j > l\) if and only if~\(W_{j+1},\ldots,W_{j+l}\) are empty, which happens with
probability~\((1-l\cdot\frac{A^2_k}{k})^{km} \leq e^{-lA_k m}\leq e^{-lm}\) since~\(A_k \geq 1.\)  Letting~\(\alpha_0\) be the smallest integer greater than or equal to~\(\alpha\) we use the moment estimate~(\ref{disc_tel}) with~\(\theta = m \geq 1\) and~\(r = 2\alpha_0\) to get~
\begin{equation}\label{t_next_est}
\mathbb{E}(T^{next}_j)^{2\alpha} \leq \mathbb{E}(T^{next}_j)^{2\alpha_0} \leq \frac{(2\alpha_0)!}{(1-e^{-m})^{2\alpha_0}} \leq \frac{(2\alpha_0)!}{(1-e^{-1})^{2\alpha_0}}.
\end{equation}

Next, if the event~\(E^c_{nice}(j)\) occurs, then in each pair~\(\{W^{small}_j(i),W^{small}_{j+1}(i)\}\)\\\( 1 \leq i \leq L\) at least one rectangle is empty. If there are~\(l \geq L\) such empty rectangles in total, the area formed by these rectangles is~\(l \frac{4A_k^2 \log{m}}{km}\) and since there are~\(2L\) rectangles to choose from, this happens with probability at most~\[{2L \choose l}\left(1-l \cdot \frac{4A_k^2 \log{m}}{km}\right)^{km} \leq (2L)^{l}e^{-4lA_k^2\log{m}} \leq (2L)^{l}e^{-4l\log{m}}\] since~\(A_k \geq 1.\) Using~\(L \leq \sqrt{m},\) we further have~\[(2L)^{l}e^{-4l\log{m}} \leq (2\sqrt{m})^{l} e^{-4l\log{m}}  \leq e^{-3l\log{m}}\] since~\(2\sqrt{m}e^{-4\log{m}} = \frac{2\sqrt{m}}{m^{4}} \leq \frac{1}{m^3}\) for all~\(m \geq 4.\)

Therefore
\[\mathbb{P}(E^c_{nice}(j)) \leq \sum_{l \geq L} e^{-3l\log{m}} \leq De^{-3L\log{m}} \leq De^{-2\sqrt{m}\log{m}}, \] for some constant~\(D > 0\) not depending on~\(j,k\) or~\(m,\) since~\(L \geq \frac{\sqrt{m}}{2}.\)
Combining this and the estimate~(\ref{t_next_est}) for~\(T^{next}_j\) from the previous paragraph,
we have from~(\ref{ed_alp2}) that~\(\mathbb{E}\left(T^{next}_j\right)^{\alpha}\ind(E^c_{dense}(j))  \leq D_3 e^{-2\sqrt{m}\log{m}}\)
for some constant~\(D_3 >0\) not depending on~\(j,k\) or~\(m.\) Plugging this into~(\ref{ed_alpha}), we get~(\ref{s_eval}).~\(\qed\)

%is true.

%for all~\(m \geq M_0\) large, not depending on~\(k.\)

\emph{Proof of~(\ref{i1_est}) in Lemma~\ref{i1_i2_lem}}: We write~\[I_1 = \sum_{j = j_{low}}^{j_{up}} \mathbb{E} MST(N(1)) \ind(N(1) = j),\]
where~\(j_{low} := mA_k^2 - \sqrt{m}\log{m} \leq mA_k^2 + \sqrt{m}\log{m} =: j_{up}.\) Given~\(N(1) = j,\) the nodes in~\(W_1\) are independently and uniformly distributed in~\(W_1\) and we define~\(\mathbb{E}MST\left(j;\frac{A_k}{\sqrt{k}}\right)\) to be the
expected length of the shifted MST of the~\(j\) nodes independently and uniformly distributed in the~\(\frac{A_k}{\sqrt{k}} \times \frac{1}{k}\) square~\(W_1.\) From the scaling relation~(\ref{mst_scale}) in Appendix we have that
\begin{equation}
I_1 = \sum_{j =j_{low}}^{j_{up}} \mathbb{E}MST\left(j;\frac{1}{k}\right) \mathbb{P}(N(1) = j) = \left(\frac{A_k}{\sqrt{k}}\right)^{\frac{\alpha}{2}} \sum_{j=j_{low}}^{j_{up}} \left(\mathbb{E}MST_j \right)\mathbb{P}(N(1) = j), \label{temp1}
\end{equation}
by~(\ref{mst_scale}).

From the one node difference estimate~(\ref{one_diff_est_mst}), we have for~\(j_{low} \leq u \leq j_{up}\) that~
\(\mathbb{E}|MST_{u+1}-MST_{u}| \leq \left(\frac{D}{u}\right)^{\frac{\alpha}{2}}\) for some constant~\(D >0\) not depending on~\(u,k\) or~\(m.\)
Since~\(u \geq j_{low} = mA_k^2 - \sqrt{m} \log{m} \geq m-\sqrt{m}\log{m}\) using~\(A_k \geq 1\) we have
for any~\(j_{low} \leq j_1, j_2 \leq j_{up}\) that~\(\mathbb{E}|MST_{j_2} - MST_{j_1}| \leq \sum_{u=j_{low}}^{j_{up}-1} \mathbb{E}|MST_{u+1} - MST_{u}| \) is bounded above by
\begin{equation}
\sum_{u=j_{low}}^{j_{up}-1} \left(\frac{D}{u}\right)^{\frac{\alpha}{2}} \leq (j_{up}-j_{low})\left(\frac{D_1}{m}\right)^{\frac{\alpha}{2}} \leq D_2(\log{m})m^{\frac{1-\alpha}{2}},\label{temp3}
\end{equation}
where~\(D_1,D_2 > 0\) are constants not depending on~\(j_1,j_2,k\) or~\(m.\) Setting~\(j_1= mA_k^2\) and~\(j_2 = j\) and using~(\ref{temp3}) we get~\(MST_j \leq MST_{mA_k^2} + D_2 m^{\frac{1-\alpha}{2}}(\log{m})\) for all~\(j_{low} \leq j \leq j_{up}.\) From the expression for~\(I_1\) in~(\ref{temp1}) we therefore have that
\begin{equation}
I_1 \leq \left(\frac{A_k}{\sqrt{k}}\right)^{\alpha} \left(\mathbb{E} MST_{mA_k^2} + D_2 m^{\frac{1-\alpha}{2}} (\log{m})\right). \label{temp6}
\end{equation}
Since~\(A_k \leq 1+\frac{1}{\log{k}},\) we have that~\(mA_k^2 - m \leq \frac{3m}{\log{k}}\) for all~\(k \geq 4\) and so using the estimate~(\ref{mst_ab2}) with~\(n_1 = m, n_2 = m(A_k^2 - 1) \leq \frac{3m}{\log{k}},\) we get that~\(\mathbb{E} MST_{mA_k^2} \leq \mathbb{E}MST_{m} + \frac{D_3m^{1-\frac{\alpha}{2}}}{\log{k}}\) for some constant~\(D_3 >0\) not depending on~\(m\) or~\(k\)
and this obtains~(\ref{i1_est}).~\(\qed\)

%since there are at most~\(2\log{m}\) terms in the summation term in the first line of~(\ref{temp6}) (see~(\ref{j_range})).

\emph{Proof of~(\ref{i2_est}) in Lemma~\ref{i1_i2_lem}}: There are~\(N(1)\) nodes in the square~\(W_1\) and given~\(N(1) = l,\) the~\(l\) nodes are uniformly distributed within~\(W_1\) and so arguing as in the previous two paragraphs and using the scaling relation~(\ref{mst_scale}) in Appendix,
we have
\begin{eqnarray}
\mathbb{E}\left(MST(N(1))|N(1) =l\right) &=& \mathbb{E}MST\left(l;\frac{A_k}{\sqrt{k}}\right) \nonumber\\
&=& \left(\frac{A_k}{\sqrt{k}}\right)^{\alpha} \mathbb{E}MST_l \nonumber\\
&\leq& D_1\left(\frac{A_k}{\sqrt{k}}\right)^{\alpha}\cdot l^{1-\frac{\alpha}{2}}
\end{eqnarray}
for some constant~\(D_1 > 0\) using the expectation upper bound in~(\ref{exp_mst_bound}).
Thus~
\begin{eqnarray}
\mathbb{E}MST(N(1))\ind(N(1) = l) &=& \mathbb{E}(MST(N(1))|N(1)=l)\mathbb{P}(N(1) = l) \nonumber\\
&\leq& D_1\left(\frac{A_k}{\sqrt{k}}\right)^{\alpha}\cdot l^{1-\frac{\alpha}{2}}\mathbb{P}(N(1)= l) \nonumber
\end{eqnarray}
and consequently~\(I_2 = \mathbb{E}MST(N(1))\ind(F_1^c)\) equals
\begin{eqnarray}
&&\sum_{l\leq mA_k^2 -\sqrt{m}\log{m}} + \sum_{l \geq mA_k^2+\sqrt{m} \log{m}} \mathbb{E}MST(N(1))\ind(N(1) = l) \nonumber\\
&&\leq\;\;\; \sum_{l\leq mA_k^2 -\sqrt{m}\log{m}} + \sum_{l \geq mA_k^2 +\sqrt{m}\log{m}} D_1\left(\frac{A_k}{\sqrt{k}}\right)^{\alpha}\cdot l^{1-\frac{\alpha}{2}}\mathbb{P}(N(1)= l) \nonumber\\
&&=\;\;\; D_1\left(\frac{A_k}{\sqrt{k}}\right)^{\alpha} \mathbb{E}(N(1))^{1-\frac{\alpha}{2}}\ind(F_1^c). \label{i2_alt_est}
\end{eqnarray}

Each node~\(X_i, 1 \leq i \leq km\) has a probability~\(\frac{A^2_k}{k}\) of being present in the~\(\frac{A_k}{\sqrt{k}} \times \frac{A_k}{\sqrt{k}}\) square~\(W_1.\) Therefore the number of nodes~\(N(1)\) in the square~\(W_1\) is binomially distributed with~
\begin{equation}\label{mean_n1_est}
m \leq \mathbb{E}N(1) = mA_k^2  \leq 4m\;\;\text{ and }\;\;var(N(1)) \leq km \cdot \frac{A_k^2}{k} = mA_k^2 \leq 4m,
\end{equation}
since~\(1 \leq A_k \leq 1+\frac{1}{\log{k}} \leq 2.\) We recall from discussion prior to~(\ref{mst_km3}) that~\(F_1^c = \{N(1) \leq mA_k^2 - \sqrt{m} \log{m}\} \cup \{N(1) \geq mA_k^2 + \sqrt{m} \log{m}\}\) and so from Chebychev's inequality we get that~\(\mathbb{P}(F^c_1)  \leq \frac{var(N(1))}{(\sqrt{m} \log{m})^2} \leq \frac{4}{(\log{m})^2}.\) We show below that
\begin{equation}\label{en1_est}
\mathbb{E}(N(1))^{1-\frac{\alpha}{2}}\ind(F_1^c) \leq D_2e^{-D_3 m} + m^{1-\frac{\alpha}{2}}\mathbb{P}(F_1^c) \leq D_2e^{-D_3 m} + m^{1-\frac{\alpha}{2}}\frac{4}{(\log{m})^2}
\end{equation}
for some constants~\(D_2,D_3\) not depending on~\(k\) or~\(m\) and substituting~(\ref{en1_est}) into~(\ref{i2_alt_est}), we
then get~(\ref{i2_est}).

To prove~(\ref{en1_est}) we use~\(1 \leq A_k \leq 2\) to get that~\[\frac{m}{2} \leq mA_k^2 - \sqrt{m} \log{m} \leq mA_k^2 + \sqrt{m} \log{m} \leq 4m+\sqrt{m}\log{m} \leq 5m\] and so~\(\ind(F_1^c) \leq \ind(H_1) + \ind(H_2) + \ind(H_3) + \ind(H_4),\) where~\[H_1 := \left\{N(1) \leq \frac{m}{2}\right\}, H_2 := \left\{\frac{m}{2} \leq N(1)  \leq mA_k^2 - \sqrt{m} \log{m}\right\},\]
\[H_3 := \left\{mA_k^2 + \sqrt{m} \log{m} \leq N(1) \leq 5m\right\} \text{ and } H_4 := \{N(1) \geq 5m\}.\]

Using the bounds~\(m \leq \mathbb{E}N(1) \leq 4m\) from~(\ref{mean_n1_est}) and the standard deviation estimates~(\ref{std_dev_up}) and~(\ref{std_dev_down}) in Appendix we have that~\(\max\left(\mathbb{P}(H_1),\mathbb{P}(H_4)\right) \leq e^{-2Dm}\) for some constant~\(D > 0\) not depending on~\(k\) or~\(m\) and so
\begin{equation}\label{h1_est}
\mathbb{E}N(1)^{1-\frac{\alpha}{2}} \ind(H_1) \leq \mathbb{E}N(1)\ind(H_1) \leq \frac{m}{2}\mathbb{P}(H_1) \leq me^{-2Dm}.
\end{equation}
Using Cauchy-Schwarz inequality and the bound~\(\mathbb{E}N^2(1) = var(N(1)) + (\mathbb{E}N(1))^2 \leq 4m + (4m)^2 \leq 20m^2\) (see~(\ref{mean_n1_est}))
we also get
\begin{equation}\label{h4_est}
\mathbb{E}N(1)^{1-\frac{\alpha}{2}} \ind(H_4) \leq \mathbb{E}N(1)\ind(H_4) \leq \left(\mathbb{E}N^2(1)\right)^{\frac{1}{2}} \left(\mathbb{P}(H_4)\right)^\frac{1}{2} \leq m\sqrt{20} e^{-Dm}.
\end{equation}
Finally, for the range of~\(N(1)\) in the events~\(H_2\) and~\(H_3\) we have
\begin{equation}\label{h23_est}
\mathbb{E}N(1)^{1-\frac{\alpha}{2}} \ind(H_2 \cup H_3) \leq D_1 m^{1-\frac{\alpha}{2}}\mathbb{P}(H_2 \cup H_3)
\end{equation}
for some constant~\(D_1 > 0\) not depending on~\(k\) or~\(m.\) Adding~(\ref{h1_est}),~(\ref{h4_est}) and~(\ref{h23_est}) and using the fact that~\(H_2 \cup H_3 \subseteq F_1^c,\) we get~(\ref{en1_est}).~\(\qed\)

%\begin{eqnarray}
%|MST_{k+1} - MST_{k}| \leq C_1 r_{k} (\log{k}) \ind(Y_{tot}(n)) + k \sqrt{2} \ind(Y_{tot}^c(n)) \nonumber\\
%&\leq& C_2 \frac{(\log{k})^{3/2}}{\sqrt{k}} \ind(Y_{tot}(n)) + (n+1)^2 \sqrt{2} \ind(Y_{tot}^c(n)) \nonumber\\
%\label{k_dif}
%\end{eqnarray}
%for each~\(n^2 \leq k <(n+1)^2\) and for some constants~\(C_1,C_2 > 0\) not depending on~\(k\) or~\(n.\)
%The final estimate~(\ref{k_dif}) is true since from~(\ref{rn_def_mst}) we get that~\(r_k \leq C_3 \sqrt{\frac{\log{k}}{k}}\)
%for some constant~\(C_3 > 0.\)

%at the beginning of Section~\ref{pf_mst}).
%Dividing the unit square into strips of size~\(\frac{1}{\sqrt{n}} \times 1\) and arguing as in the proof of~\((b4),\)
%we obtain~\(MST_n \leq \sqrt{n} + n\frac{1}{\sqrt{n}} \sqrt{2} + 2 \leq 3\sqrt{n}\) for all~\(n\) large.

%\otimes_{k=1}^{\infty} \Omega_k

%\emph{Proof of~\((a1)-(a2)\)}:
%WRTE MORE HERE +eTC...

%We now construct an approximation~\({\cal U}_n\) of~\({\cal P}_n\) and determine the number of edges in~\({\cal U}_n.\) We write...3332

\setcounter{equation}{0}
\renewcommand\theequation{\thesection.\arabic{equation}}
\section*{Appendix : Miscellaneous results}\label{appendix}
\subsection*{\em Standard deviation estimates}
We use the following standard deviation estimates for sums of independent Poisson and Bernoulli random variables (see Alon and Spencer (2008)).
\begin{Lemma}\label{app_lem}
Suppose~\(W_i, 1 \leq i \leq m\) are independent Bernoulli random variables satisfying~\(\mu_1 \leq \mathbb{P}(W_1=1) = 1-\mathbb{P}(W_1~=~0) \leq \mu_2.\) For any~\(0 < \epsilon < \frac{1}{2},\)
\begin{equation}\label{std_dev_up}
\mathbb{P}\left(\sum_{i=1}^{m} W_i > m\mu_2(1+\epsilon) \right) \leq \exp\left(-\frac{\epsilon^2}{4}m\mu_2\right)
\end{equation}
and
\begin{equation}\label{std_dev_down}
\mathbb{P}\left(\sum_{i=1}^{m} W_i < m\mu_1(1-\epsilon) \right) \leq \exp\left(-\frac{\epsilon^2}{4}m\mu_1\right)
\end{equation}
Estimates~(\ref{std_dev_up}) and~(\ref{std_dev_down}) also hold if~\(\{W_i\}\) are independent Poisson random variables with~\(\mu_1 \leq \mathbb{E}W_1 \leq \mu_2.\)
\end{Lemma}

\subsection*{\em Proof of the monotonicity property~(\ref{mon_salpha})}
For~\(\alpha \leq 1\) we couple the original Poisson process~\({\cal P}\) and the homogenous process~\({\cal P}_{\delta}\) in the following way. Let~\(V_{i}, i \geq 1\) be i.i.d.\ random variables each with density~\(f(.)\) and let~\(N_V\) be a Poisson random variable with mean~\(n,\) independent of~\(\{V_i\}.\) The nodes~\(\{V_i\}_{1 \leq i \leq N_V}\) form a Poisson process with intensity~\(nf(.)\) which we denote as~\({\cal P}\) and colour green.

Let~\(U_i, i \geq 1\) be i.i.d.\ random variables each with density~\(\epsilon_2-f(.)\) where~\(\epsilon_2 \geq 1\) is as in~(\ref{f_eq}) and let~\(N_U\) be a Poisson random variable with mean~\(n(\epsilon_2-1).\) The random variables~\((\{U_i\},N_U)\) are independent of~\((\{V_i\},N_V)\) and the nodes~\(\{U_i\}_{1 \leq i  \leq N_U}\) form a Poisson process with intensity~\(n(\epsilon_2-f(.))\) which we denote as~\({\cal P}_{ext}\) and colour red. The nodes of~\({\cal P}\) and~\({\cal P}_{ext}\) together form a homogenous Poisson process with intensity~\(n\epsilon_2,\) which we denote as~\({\cal P}_{\delta}\) and define it on the probability space~\((\Omega_{\delta},{\cal F}_{\delta}, \mathbb{P}_{\delta}).\)

Let~\(\omega_{\delta} \in \Omega_{\delta}\) be any configuration and as above let~\(\{i^{(\delta)}_{j}\}_{1 \leq j \leq Q_{\delta}}\) be the indices of the squares in~\(\{R_j\}\) containing at least one node of~\({\cal P}_{\delta}\) and let~\(\{i_{j}\}_{1 \leq j \leq Q}\) be the indices of the squares in~\(\{R_j\}\) containing at least one node of~\({\cal P}.\) The indices in~\(\{i^{(\delta)}_j\}\) and~\(\{i_j\}\) depend on~\(\omega_{\delta}.\) Defining~\(S_{\alpha} = S_{\alpha}(\omega_{\delta})\) and~\(S^{(\delta)}_{\alpha} = S_{\alpha}^{(\delta)}(\omega_{\delta})\) as before, we have that~\(S_{\alpha}\) is determined only by the green nodes of~\(\omega_{\delta}\) while~\(S^{(\delta)}_{\alpha}\) is determined by both green and red nodes of~\(\omega_{\delta}.\)

From the monotonicity property, we therefore have that~\(S_{\alpha}(\omega_{\delta}) \leq S^{(\delta)}_{\alpha}(\omega_{\delta})\) and so for any~\(x > 0\) we have
\begin{equation}\label{mon_eq}
\mathbb{P}_{\delta}(S^{(\delta)}_{\alpha} < x) \leq \mathbb{P}_{\delta}(S_{\alpha} < x)  = \mathbb{P}_0(S_{\alpha} < x),
\end{equation}
proving~(\ref{mon_salpha}).

If~\(\alpha > 1,\) we perform a slightly different analysis. Letting~\(\epsilon_1 \leq 1\) be as in~(\ref{f_eq}), we construct a Poisson process~\({\cal P}_{ext}\) with intensity~\(n(f(.)-\epsilon_1)\) and colour nodes of~\({\cal P}_{ext}\) red. Letting~\({\cal P}_{\delta}\) be another independent Poisson process with intensity~\(n\epsilon_1,\) we colour nodes of~\({\cal P}_{\delta}\) green. The superposition of~\({\cal P}_{ext}\) and~\({\cal P}_{\delta}\) is a Poisson process with intensity~\(nf(.),\) which we define on the probability space~\((\Omega_{\delta},{\cal F}_{\delta},\mathbb{P}_{\delta}).\) In this case, the sum~\(S_{\alpha}\) is determined by both green and red nodes while~\(S^{(\delta)}_{\alpha}\) is determined only by the green nodes. Again using the monotonicity property of~\(S_{\alpha},\) we get~(\ref{mon_eq}).~\(\qed\)

\subsection*{\em Scaling and translation property of MSTs}
For a set of nodes~\(\{x_1,\ldots,x_n\}\) in the unit square~\(S,\)
recall from Section~\ref{intro} that~\(K_n(x_1,\ldots,x_n)\)
is the complete graph formed by joining all the nodes by straight line segments and the edge~\((x_i,x_j)\)
is assigned a weight of~\(h^{\alpha}(x_i,x_j),\)
where~\(h(x_i,x_j)\) is the weight of the edge~\((x_i,x_j)\) satisfying~(\ref{eq_met1}) and properties~\((b1)-(b2).\)
We denote~\(MST(x_1,\ldots,x_n)\) to be the length of the minimum spanning tree of~\(K_n(x_1,\ldots,x_n)\)
with edge weights obtained as in~(\ref{min_weight_tree}).

\emph{Scaling}: For any~\(a > 0,\) consider the graph~\(K_n(ax_1,\ldots,ax_n).\) Using the scaling property~\((b1),\) the weight~\(h(ax_1,ax_2)\) of the edge between
the vertices~\(ax_1\) and~\(ax_2\) is simply~\(a \cdot h(x_1,x_2),\) where~\(h(x_1,x_2)\) is the weight of the edge between~\(x_1\) and~\(x_2.\)
Using the definition of MST in~(\ref{min_weight_tree})
we therefore have~\(MST(ax_1,\ldots,ax_n) = a^{\alpha} MST(x_1,\ldots,x_n)\)
and so if~\(Y_1,\ldots,Y_n\) are~\(n\) nodes uniformly distributed in the square~\(aS\) of side length~\(a,\)
then~\[MST(n;a) := MST(Y_1,\ldots,Y_n) = a^{\alpha} MST(X_1,\ldots,X_n),\]
where~\(X_i = \frac{Y_i}{a}, 1 \leq i \leq n\) are i.i.d.\ uniformly distributed in~\(S.\)
Recalling the notation~\(MST_n = MST(X_1,\ldots,X_n)\) from~(\ref{min_weight_tree}) we therefore get
\begin{equation}\label{mst_scale}
\mathbb{E} MST(n;a) = a^{\alpha} \mathbb{E}MST_n.
\end{equation}

\emph{Translation}: For~\(b \in \mathbb{R}^2\) consider the graph~\(K_n(x_1+b,\ldots,x_n+b).\) Using the translation property~\((b2),\)
the weight~\(h(x_1+b,x_2+b) \leq h_0 \cdot h(x_1,x_2),\) the weight of the edge between~\(x_1\) and~\(x_2.\)
Using the definition of MST in~(\ref{min_weight_tree})
we therefore have~\(MST(x_1+b,\ldots,x_n+b) \leq h_0^{\alpha} \cdot MST(x_1,\ldots,x_n).\)~\(\qed\)

\subsection*{\em Moments of random variables}
Let~\(X \geq 1\) be any integer valued random variable such that
\begin{equation}\label{x_dist}
\mathbb{P}(X \geq l) \leq e^{-\theta (l-1)}
\end{equation}
for all integers~\(l \geq 1\) and some constant~\(\theta > 0\) not depending on~\(l.\) For every integer~\(r \geq 1,\)
\begin{equation}\label{disc_tel}
\mathbb{E}X^{r} \leq r\sum_{l\geq 1}  l^{r-1} \mathbb{P}(X \geq l) \leq r\sum_{l \geq 1} l^{r-1} e^{-\theta (l-1)} \leq \frac{r!}{(1-e^{-\theta})^{r}}
\end{equation}
\emph{Proof of~(\ref{disc_tel})}: For~\(r \geq 1\) we have
\begin{equation}\label{tel_sum}
\mathbb{E}X^{r} = \sum_{l \geq 1} l^{r} \mathbb{P}(X= l) = \sum_{l \geq 1}l^{r} \mathbb{P}(X \geq l) - l^{r}\mathbb{P}(X \geq l+1)
\end{equation}
and substituting the~\(l^{r}\) in the final term of~(\ref{tel_sum}) with~\(  (l+1)^{r} - ((l+1)^{r}-l^{r})\) we get
\begin{eqnarray}
\mathbb{E}X^r &=& \sum_{l \geq 1} \left(l^{r} \mathbb{P}(X \geq l) - (l+1)^{r} \mathbb{P}(X \geq l+1)\right)  \nonumber\\
&&\;\;\;\;\;\;\;\; + \;\;\;\sum_{l \geq 1} ((l+1)^{r}-l^{r}) \mathbb{P}(X \geq l+1) \nonumber\\
&=& 1 + \sum_{l \geq 1}((l+1)^{r}-l^{r}) \mathbb{P}(X \geq l+1)  \nonumber\\
&=& \sum_{l \geq 0}((l+1)^{r}-l^{r}) \mathbb{P}(X \geq l+1)  \label{gent}
\end{eqnarray}
where the second equality is true since~\(l^{r}\mathbb{P}(X \geq l) \leq l^{r}e^{-\theta (l-1)} \longrightarrow 0\) as~\(l~\rightarrow~\infty.\) Using~\((l+1)^{r} - l^{r} \leq r\cdot (l+1)^{r-1}\) in~(\ref{gent}), we get the first relation in~(\ref{disc_tel}).

We prove the second relation in~(\ref{disc_tel}) by induction as follows. Let~\(\gamma = e^{-\theta} < 1\) and~\(J_r := \sum_{l \geq 1} l^{r-1} \gamma^{l-1}\) so that
\begin{equation}
J_{r+1}(1 - \gamma) = \sum_{l \geq 1} l^{r} \gamma^{l-1} - \sum_{l \geq 1}l^{r} \gamma^{l} = \sum_{l \geq 1} \left(l^{r}-(l-1)^{r} \right)\gamma^{l-1}. \nonumber
\end{equation}
Using~\(l^{r}-(l-1)^r \leq r\cdot l^{r-1}\) for~\(l \geq 1\) we therefore get that
\[J_{r+1}(1-\gamma) \leq r\sum_{l \geq 1}l^{r-1} \gamma^{l-1} = rJ_r\]
and so the second relation in~(\ref{disc_tel}) follows from induction.~\(\qed\)

\subsection*{\em Acknowledgement}
I thank Professors Rahul Roy, C. R. Subramanian and Federico Camia for crucial comments that led to an improvement
of the paper. I also thank Professors Rahul Roy, C. R. Subramanian, Federico Camia and IMSc for my fellowships.

\bibliographystyle{plain}

\end{document}